\documentclass[12pt]{amsart}
\usepackage{amsmath,amscd,amssymb,amsfonts}
\newcommand{\h}{\hbox}
\newcommand{\q}{\quad}
\newcommand{\nin}{\noindent}
\newcommand{\bs}{\par\bigskip}
\newcommand{\ms}{\par\medskip}
\newcommand{\sk}{\par\smallskip}
\newcommand{\bsn}{\par\bigskip\noindent}
\newcommand{\msn}{\par\medskip\noindent}
\newcommand{\skn}{\par\smallskip\noindent}
\newcommand{\ssb}{\raise.15ex\h{${\scriptscriptstyle\bullet}$}}
\newcommand{\ssc}{\,\raise.15ex\h{${\scriptstyle\circ}$}\,}
\newcommand{\msum}{\hbox{$\sum$}}
\newcommand{\mprod}{\hbox{$\prod$}}
\newcommand{\mcup}{\hbox{$\bigcup$}}
\newcommand{\mopl}{\hbox{$\bigoplus$}}
\newcommand{\C}{{\mathbb C}}
\newcommand{\D}{{\mathbb D}}
\newcommand{\N}{{\mathbb N}}
\newcommand{\PP}{{\mathbb P}}
\newcommand{\Q}{{\mathbb Q}}
\newcommand{\R}{{\mathbf R}}
\newcommand{\Z}{{\mathbb Z}}
\newcommand{\A}{{\mathcal A}}
\newcommand{\DD}{{\mathcal D}}

\newcommand{\Nc}{{\mathcal N}}
\newcommand{\I}{{\mathcal I}}

\newcommand{\M}{{\mathcal M}}
\newcommand{\Oc}{{\mathcal O}}
\newcommand{\Sc}{{\mathcal S}}
\newcommand{\Et}{\widetilde{E}}

\newcommand{\Xt}{\widetilde{X}}
\newcommand{\Yt}{\widetilde{Y}}
\newcommand{\x}{\widetilde{x}}
\newcommand{\y}{\widetilde{y}}
\newcommand{\rhot}{\widetilde{\rho}}
\newcommand{\aaa}{{\mathbf a}}
\newcommand{\bb}{{\mathbf b}}
\newcommand{\cc}{{\mathbf c}}
\newcommand{\ee}{{\mathbf e}}
\newcommand{\al}{\alpha}
\newcommand{\la}{\lambda}
\newcommand{\La}{\Lambda}
\newcommand{\lah}{\widehat{\Lambda}}
\newcommand{\dd}{\partial}
\newcommand{\ddd}{{\rm d}}
\newcommand{\dg}{{\rm d}g}
\newcommand{\dx}{{\rm d}x}
\newcommand{\dy}{{\rm d}y}
\newcommand{\da}{{\rm d}^{\alpha}}
\newcommand{\kod}{\tfrac{k}{d}}
\newcommand{\om}{\omega}
\newcommand{\oma}{\omega^{\alpha}}
\newcommand{\Om}{\Omega}
\newcommand{\lH}{{\rm log}_H}
\newcommand{\sw}{\,{\wedge}\,}
\newcommand{\Sf}{{\mathfrak S}}
\newcommand{\Gr}{{\rm Gr}}
\newcommand{\mult}{{\rm mult}}
\newcommand{\prd}{{\rm prod}}
\newcommand{\Supp}{{\rm Supp}}
\newcommand{\dett}{{\rm det}}
\newcommand{\mmu}{\rlap{$\mu$}\hskip1.1pt\h{$\mu$}}
\newcommand{\bl}{\bigl}
\newcommand{\br}{\bigr}
\newcommand{\into}{\hookrightarrow}
\newcommand{\simto}{\buildrel\sim\over\longrightarrow}
\newcommand{\onto}{\mathop{\rlap{$\to$}\hskip2pt\h{$\to$}}}
\newcommand{\ges}{\geqslant}
\newcommand{\les}{\leqslant}
\newcommand{\1}{\hskip1pt}
\begin{document}
\title[Local systems and Aomoto complexes]{Rank one local systems on complements\\of hyperplanes and Aomoto complexes}
\author{Morihiko Saito}
\address{RIMS Kyoto University, Kyoto 606-8502 Japan}
\begin{abstract} We show that the cohomology of a rank 1 local system on the complement of a projective hyperplane arrangement can be calculated by the Aomoto complex in certain cases even if the condition on the sum of the residues of connection due to Esnault et al is not satisfied. For this we have to study the localization of Hodge-logarithmic differential forms which are defined by using an embedded resolution of singularities. As an application we can compute certain monodromy eigenspaces of the first Milnor cohomology group of the defining polynomial of the reflection hyperplane arrangement of type $G_{31}$ without using a computer.
\end{abstract}
\maketitle
\centerline{\bf Introduction}
\bsn
Let $X$ be a projective hyperplane arrangement of degree $d$ in $Y:=\PP^{n-1}$, and $L$ be a nontrivial rank 1 local system on the complement $U:=Y\setminus X$. Let $\la_k$ be the eigenvalue of the local monodromy of $L$ around a hyperplane $X_k$ in $X$ ($k\in[1,d]$), see also Remark~(iii) after Theorem~(2.2) below. Choosing complex numbers $\al_k$ so that
$$\exp(-2\pi i\al_k)=\la_k,\q\msum_{k=1}^d\1\al_k=0,
\leqno(1)$$
we can get a free $\Oc_U$-module of rank 1 with an integrable connection $\nabla$ calculating the local system $L$. (Note that $\prod_{k=1}^d\la_k=1$ by restricting to a generic line.) Combined with a well-known assertion about the cohomology ring $\A^{\ssb}=H^{\ssb}(U,\C)$ (see \cite[Lemma 5]{Br} and also \cite{OS}, \cite{OT}), this gives the {\it Aomoto complex} $(\A^{\ssb},\oma\wedge)$ (where $\oma\in\A^1$ is determined by $\al=(\al_k)$, see (2.1) below) together with the canonical morphisms
$$H^j(\A^{\ssb},\oma\wedge)\to H^j(U,L)\q\q(j\in\N).
\leqno(2)$$
These are isomorphisms if a certain condition for the $\al_k$ is satisfied as in \cite{ESV}; more precisely, if the following condition due to \cite{STV} holds (see Theorem~(2.2) below):
$$\al_Z:=\msum_{k\in I_Z}\,\al_k\notin\Z_{>0}\,\,\,\h{for any {\it dense edge} $Z$,}
\leqno(3)$$
where $I_Z:=\{k\in[1,d]\mid X_k\supset Z\}$. (Note that $\al\ne0$, since $L$ is assumed nontrivial.)
\sk
The above isomorphisms are quite useful to calculate the Milnor cohomology and the Bernstein-Sato polynomial of the defining polynomial of $X$ if we can take the $\al_k$ so that conditions (1) and (3) are satisfied, see (2.5) below. However, these conditions are fairly strong, and it cannot be done in general unless the arrangement is rather simple, see for instance \cite{BDS}, \cite{BSY}, \cite{bCM}, \cite{bha}. An idea to relax conditions (1) and (3) was suggested in the proof of \cite[2.2.1]{bha}. It seems, however, rather difficult to realize it {\it in a simple way}. We may do it, for instance, as follows.
\sk
Assume $n=3$. If we are interested in the calculation of the {\it first} Milnor cohomology, we may assume it as is well-known using (2.5) below and a hyperplane section cut together with the Artin (or weak Lefschetz type) vanishing theorem in \cite{BBD}, \cite{Di2}, see also \cite[Lemma 2.1.18]{mhp}. In the case $n=3$, each dense edge $Z$ of codimension 2 in $Y=\PP^2$ is identified with a singular point $p$ of $X$ with multiplicity at least 3 (where $X$ is assumed reduced). Set
$$\al_p:=\al_Z,\q I_p:=I_Z\q\h{if}\q Z=\{p\}.$$
Let $p_i$ ($i\in[1,r]$) be the points of $X$ such that condition~(3) is {\it not\1} satisfied for $Z=\{p_i\}$. Set
$$\aligned&I^b:=\mcup_{i=1}^r\,I_{p_i},\q I^c:=\{1,\dots,d\}\setminus I^b,\\&X^{I^b}:=\mcup_{k\in I^b}\,X_k,\q X^{I^c}:=\mcup_{k\in I^c}\,X_k.\endaligned$$
Here $I^b$ means {\it bad\1} indices. Let $X^{(\beta)}$ be the set of singular points $p$ of $X$ with multiplicity at least 3 and such that $\al_p=\beta$, where $\beta$ is any complex number. We define similarly $X^{I^c(\beta)}$ by replacing $X$ with $X^{I^c}$. We have the following.
\msn
{\bf Theorem 1.} {\it Assume $n=3$. The morphism $(2)$ is bijective for $j=1$ and the equality of the dimensions of the source and target of $(2)$ holds for any $j$, if the following conditions are satisfied\,$:$
\skn
{\rm (a)} $p_i\in X_d$, $\,\al_{p_i}=1\,\,\,(\forall\,i\in[1,r]),\q\al_d\ne 1$.
\skn
{\rm (b)} $X^{I^c}\setminus\bl(X^{I^c(0)}\cup X^{I^b}\br)$ is connected.
\skn
{\rm (c)} $\#\bl(X_k\cap X^{I^c}\setminus\mcup_{j\in I^b\setminus\{k\}}X_j\br)\ges 2\,\,\,\,(\forall\,k\in I^b)$,\q$\al_j\ne 0\,\,\,\,(\forall\,j\in I^c)$.}
\skn
{\rm (d)} $X_d$ contains at least two ordinary double points of $X$.
\ms
Note that condition~(c) can be weakened slightly, see Remark after (3.4) below.
The equalities of dimensions follows from the relation between the Euler characteristic of $U$ and that of each monodromy eigenspace of the Milnor cohomology groups, see (2.5.1) below. As an application, we can prove the vanishing of the monodromy eigenspace of the first Milnor cohomology $H^1(F_{\!f},\C)_{\la}$ with $\la=\exp(-2\pi i/6)$ for the reflection hyperplane arrangement of type $G_{31}$ as in (A.5) below, see \cite{BDY}, \cite{DS} for proofs using a computer, \cite{MPP} for other eigenspaces (that is, for $\la$ satisfying $\la^p=1$ with $p$ prime, where $p=2$ or $3$ in this case), and \cite{Di3} for other reflection arrangements. We can show the assertion also for $\la=-1$ using Theorem~1, but not for $\la=\exp(\pm 2\pi i/3)$, see Remarks (ii) and (iii) after (A.5). If $d$ is not sufficiently large as in the case of $G_{31}$ (with $d=60$), it is not necessarily easy to satisfy conditions~(b) and (c), and some modification may be needed, see Example~(2.6)(i) below.
\sk
This work is partially supported by Kakenhi 15K04816. I would like to thank A.~Dimca for drawing my attention to this subject.
\sk
In Section~1 we review Hodge-logarithmic differential forms along hyperplane arrangements in projective spaces, which are defined by using an embedded resolution. In Section~2 we review Aomoto complexes which can calculate the cohomology of rank 1 local systems on the complements of hyperplane arrangements under some hypothesis. In Section~3 we restrict to the case $n=3$, and prove Theorem~1 for the calculation of the cohomology of a rank 1 local system via the Aomoto complex. In Appendix we study some combinatorics of the intersection poset of the reflection hyperplane arrangement of type $G_{31}$.
\bs\bs
\vbox{\centerline{\bf 1. Hodge-logarithmic differential forms}
\bsn
In this section we review Hodge-logarithmic differential forms along hyperplane arrangements in projective spaces, which are defined by using an embedded resolution.}
\msn
{\bf 1.1.~Construction.} Let $X\subset Y:=\PP^{n-1}$ be a projective hyperplane arrangement. We denote the set of edges of $X$ with dimension $j$ and codimension $n-1-j$ by
$$\Sc_{X,j}=\Sc_X^{n-1-j}\q\q(j\in[0,n{-}1]).$$
Here an {\it edge} means an intersection of hyperplanes in $X$, and $Y$ itself is viewed as an edge with $j=n{-}1$. Set
$$\Sc_X:=\h{$\bigsqcup$}_{j=0}^{n-1}\,\Sc_{X,j}=\h{$\bigsqcup$}_{j=0}^{n-1}\,\Sc_X^{n-1-j}.$$
We have a sequence of blow-ups
$$\rho:\Yt=Y_{n-2}\buildrel{\rho_{n-3}}\over\longrightarrow\,\,\cdots\,\,\buildrel{\rho_{j+1}}\over\longrightarrow Y_{j+1}\buildrel{\rho_j}\over\longrightarrow Y_j\longrightarrow\,\,\cdots\,\,\buildrel{\rho_0}\over\longrightarrow Y_0=Y,$$
where $\rho_j:Y_{j+1}\to Y_j$ is the blow-up along the center $C_j\subset Y_j$ which is the {\it disjoint union} of the proper transforms of $Z\in\Sc_{X,j}$, where we do {\it not\1} restrict to dense edges as in \cite{STV}, see Remark after Proposition~(1.2) below. Set
$$\rhot_j:=\rho_j\ssc\cdots\ssc\rho_{n-3}:\Yt\to Y_j\q(j\in[0,n-3]),$$
so that
$$\rhot_j=\rho_j\ssc\rhot_{j+1}\q(j\in[0,n-4]),\q\rhot_0=\rho.$$
Let $\Xt\subset\Yt$ be the total transforms of $X$. This is a divisor with simple normal crossings.
\msn
{\bf Proposition~1.2.} {\it We have}
$$R^k\rho_*\Om_{\Yt}^p(\log\Xt)=0\q(k>0,\,p\in\Z).
\leqno(1.2.1)$$
\msn
{\it Proof.} We prove the following by decreasing induction on $j\in[0,n-2]$:
$$R^k(\rhot_j)_*\Om_{\Yt}^p(\log\Xt)=0\q(k>0,\,p\in\Z).
\leqno(1.2.2)$$
We first show
$$\Supp\,R^k(\rhot_j)_*\Om_{\Yt}^p(\log\Xt)\subset C_j\q(k>0).
\leqno(1.2.3)$$
This is trivial if $j=n-2$. For $j\in[0,n-3]$, we have the spectral sequence
$$E_2^{k,i}=R^k(\rho_j)_*R^i(\rhot_{j+1})_*\Om_{\Yt}^p(\log\Xt)\Longrightarrow R^{k+i}(\rhot_j)_*\Om_{\Yt}^p(\log\Xt),
\leqno(1.2.4)$$
degenerating at $E_2$, since we have by inductive hypothesis
$$R^i(\rhot_{j+1})_*\Om_{\Yt}^p(\log\Xt)=0\q(i>0).$$
We then get the canonical isomorphisms
$$R^k(\rhot_j)_*\Om_{\Yt}^p(\log\Xt)=R^k(\rho_j)_*(\rhot_{j+1})_*\Om_{\Yt}^p(\log\Xt)\q(k\ges 0).
\leqno(1.2.5)$$
So (1.2.3) follows.
\sk
By \cite[Theorem 4.1.5]{Gr}, the completion of $R^k(\rhot_j)_*\Om_{\Yt}^p(\log\Xt)$ along $C_j$ can be calculated by the projective limit over $i\in\N$ of
$$R^k(\rhot_j)_*\bl(\Om_{\Yt}^p(\log\Xt)\otimes_{\Oc_{\Yt}}\Oc_{\Yt}/\I_{\Et_j}^{\,i}\br),$$
where $\I_{\Et_j}$ is the ideal sheaf of the proper transform $\Et_j$ of the exceptional divisor $E_j$ of the blow-up $\rho_j$. Note that $\Et_j$ coincides with the total transform of $E_j$, since the intersection of $E_j$ with the proper transform of $Z\in\Sc_{X,j'}$ ($j'>j$) is transversal (if it is non-empty).
\sk
We have the product structure of $\Et_j$ as in \cite[Proposition 2.6]{BS}. (Note that the projectified normal bundle of $C_j\subset Y_j$ is trivial.) We may then assume that the center of the blow-up is a point (considering the hyperplanes arrangement defined by hyperplanes in $X$ containing a given $Z\in\Sc_{X,j}$). The situation is thus reduced essentially to the case $j=0$, where $E_j$ is projective space.
\sk
We now omit the index $j$ to simplify the notation. Let $\Xt'$ be the closure of $\Xt\setminus\Et$. Set $\Xt'_{\Et}:=\Xt'\cap\Et$. Then $(\Et,\Xt'_{\Et})$ is an embedded resolution of the intersection of the proper transform of $X$ with $E$. By the above argument, the assertion (1.2.2) has been reduced to
$$H^k(\Et,\Om_{\Yt}^p(\log\Xt)\otimes_{\Oc_{\Yt}}\I_{\Et}^i/\I_{\Et}^{i+1})=0\q(k>0,\,i\ges 0).
\leqno(1.2.6)$$
\sk
Using local coordinates associated with the normal crossing divisor, we have the surjection
$$\Om_{\Yt}^p(\log\Xt)|_{\Et}\onto\Om_{\Yt}^p(\log\Xt)/\Om_{\Yt}^p(\log\Xt')=\Om_{\Et}^{p-1}(\log\Xt'_{\Et}).
\leqno(1.2.7)$$
Here we use analytic sheaves (together with GAGA if necessary). The last isomorphism is induced by taking the residue along $\Et$, and $|_{\Et}$ means the restriction as $\Oc$-modules (that is, the tensor product with $\Oc_{\Et}$).
\sk
The above surjection induces the short exact sequence
$$0\to\Om_{\Et}^p(\log\Xt'_{\Et})\to\Om_{\Yt}^p(\log\Xt)|_{\Et}\to\Om_{\Et}^{p-1}(\log\Xt'_{\Et})\to 0,
\leqno(1.2.8)$$
where the first morphism is constructed by using the surjection
$$\Om_{\Yt}^p(\log\Xt')|_{\Et}\onto\Om_{\Et}^p(\log\Xt_{\Et}).$$
The exactness of (1.2.8) is then verified by using local coordinates, see also \cite[4.4.3]{BS} for the assertion in the Grothendieck group.
\sk
Since $E$ is the exceptional divisor of a point center blow-up, and intersects transversally the proper transforms of edges passing through the center, we have the isomorphisms
$$\I_{\Et}^{\,i}/\I_{\Et}^{\,i+1}=\Oc_{\Et}(i)\q(i\ges 0),
\leqno(1.2.9)$$
where $(i)$ denotes the tensor product with the pull-back of $\Oc_E(i)$. Note that $E$ is projective space, and the morphism $\Et\to E$ is the embedded resolution of the intersection of $E$ with the proper transform of $X$ as in the construction in (1.1). (This is closely related to the product structure mentioned above.) So the assertion (1.2.6) follows from \cite[2.2.1]{bha} using (1.2.8--9). This finishes the proof of Proposition~(1.2).
\msn
{\bf Remark.} If we blow-up only for dense edges as in \cite{STV}, then the proper transforms of dense edges of dimension $j$ in $Y_j$ are not necessarily disjoint. We have to factorize further $\rho_j:Y_{j+1}\to Y_j$ in this case, and some more nontrivial argument is required at the intersections of the proper transforms of dense edges (that is, the compatibility of blow-ups with the product with smooth varieties).
\msn
{\bf 1.3.~Hodge-logarithmic forms.} In the notation of (1.1), we define the sheaf of {\it Hodge-logarithmic differential forms} on $Y$ by
$$\Om_Y^p(\lH X):=\rho_*\Om_{\Yt}^p(\log\Xt).
\leqno(1.3.1)$$
This sheaf can be either analytic or algebraic as the reader prefers, since there is essentially no difference between them by GAGA as long as the cohomology groups over $Y$ are considered.
\msn
{\bf Remark.} The above sheaf is different from the usual logarithmic differential forms (which are defined by the condition that $h\om$ and $h\ddd\om$ are holomorphic with $h$ a local defining function of the divisor, see \cite{SaK}). More precisely, we have the inclusion
$$\Om_Y^p(\lH X)\subset\Om_Y^p(\log X),$$
but the equality does not necessarily hold in general. For instance, consider the case
$$X=\{xy(x+y)z=0\}\subset Y=\PP^2,$$
where the sheaf $\Om_Y^2(\lH X)$ is not locally free at $0:=[0:0:1]\in\PP^2$, although $\Om_Y^2(\log X)$ is locally free. More precisely, we have
$$\Om_Y^2(\lH X)\cong\I_0(1),\q\Om_Y^2(\log X)\cong\Oc_Y(1),$$
with $\I_0\subset\Oc_Y$ the ideal sheaf of $0\in Y\,(=\PP^2)$, and
$$\dim H^0\bl(Y,\Om_Y^2(\lH X)\br)=2,\q\dim H^0\bl(Y,\Om_Y^2(\log X)\br)=3.$$
\sk
Note also that global logarithmic differential forms on $Y$ are not necessarily closed in general. For instance, in the above example, set
$$\x:=x/z,\q\y:=y/z,\q h:=\x\y(\x+\y),$$
$$\om:=(\y/h)\ddd \x-(\x/h)\ddd \y.$$
Then $\om$ is a non-closed logarithmic differential form on $Y$ along $X$, and $\ddd\om$ is a logarithmic form on $Y$ which is not Hodge-logarithmic, see also \cite[Remark A.3\,(ii)]{ex}, \cite{Wo}.
\msn
{\bf Corollary~1.4.} {\it We have
$$H^k\bl(Y,\Om_Y^p(\lH X)(i)\br)=0\q(k>0,\,i,p\ges 0),
\leqno(1.4.1)$$
where $(i)$ means the tensor product with $\Oc_Y(i)$ over $\Oc_Y$.}
\msn
{\it Proof.} This follows from Proposition~(1.2) and \cite[2.2.1]{bha} using the spectral sequence for the composition of cohomological functors $H^{\ssb}(Y,\,)$ and $R^{\ssb}\rho_*$ (similar to (1.2.4)).
\msn
{\bf Corollary~1.5.} {\it There are canonical isomorphisms for $p\ges 0:$}
$$\aligned H^0\bl(Y,\Om_Y^p(\lH X)\br)&=H^0\bl(\Yt,\Om_{\Yt}^p(\log\Xt)\br)=\Gr_F^pH^p(U,\C)\\&=\Gr_F^p\Gr_{2p}^WH^p(U,\C)=H^p(U,\C).\endaligned
\leqno(1.5.1)$$
\msn
{\it Proof.} The first isomorphism follows from Corollary~1.4, and the others from \cite[Lemma 5]{Br} and \cite{De}, see also \cite{ESV}.
\msn
{\bf Remarks.} (i) The last isomorphisms in (1.5.1) mean that the cohomology $H^p(U,\C)$ consists of integral logarithmic forms of type $(p,p)$ (where integral means that its cohomology class comes from the cohomology with $\Z$-coefficients), see (2.4) below for a more precise assertion.
\ms
(ii) For $p=0$, we have
$$\Om_Y^0(\lH X)=\Oc_Y,
\leqno(1.5.2)$$
since $\Om_{\Yt}^0(\log\Xt)=\Oc_{\Yt}$ and $Y,\Yt$ are smooth.
\bs\bs
\vbox{\centerline{\bf 2. Aomoto complexes}
\bsn
In this section we review Aomoto complexes which can calculate the cohomology of rank 1 local systems on the complements of hyperplane arrangements under some hypothesis.}
\msn
{\bf 2.1.~.Cohomology of rank $1$ local systems} Let $X$ be a (reduced) hyperplane arrangement in $Y=\PP^{n-1}$. By \cite[Lemma 5]{Br} and \cite{ESV}, \cite{STV}, the cohomology of a rank 1 local system $L$ on $U:=Y\setminus X$ can be calculated in certain cases as follows.
\sk
Let $X_k\,(k\in[1,d])$ be the hyperplanes in $X$ with $d:=\deg X$. 
Let $\la_k$ be the eigenvalue of the monodromy of $L$ around $X_k$. Restricting to a generic line on $Y$, we get the relation
$$\mprod_{j=1}^d\,\la_k=1.
\leqno(2.1.1)$$
\sk
Set
$$Y':=Y\setminus X_d\,\,(=\C^{n-1}),\,\,\,\,X':=X\setminus X_d,\,\,\,\,X'_k:=X_k\cap Y'\,\,\,(k\in[1,d-1]).$$
Let $g_k$ be a polynomial of degree $1$ on $Y'\,(=\C^{n-1})$ defining $X'_k$. Put
$$\om_k=\ddd g_k/g_k\q(k\in[1,d-1]).$$
Let $\al=(\al_1,\dots,\al_{d-1})\in\C^{d-1}$ satisfying
$$\exp(-2\pi i\al_k)=\la_k\q(k\in[1,d-1]).
\leqno(2.1.2)$$
Set
$$g:=\mprod_{k=1}^{d-1}\,g_k,\q g^{\al}:=\mprod_{k=1}^{d-1}\,g_k^{\al_k},\q\oma:=\msum_{k=1}^{d-1}\,\al_k\om_k.$$
Let $\Oc_{Y'}g^{\al}$ be the free $\Oc_{Y'}$-module of rank $1$ on $Y'$ generated by $g^{\al}$. Since
$$\ddd g^{\al}=g^{\al}\oma,$$
there is a regular singular integrable connection $\nabla$ on $\Oc_{Y'}g^{\al}$ satisfying
$$\nabla(hg^{\al})=g^{\al}(\ddd h)+h\1g^{\al}\oma\q\h{for}\q h\in\Oc_{Y'}.$$
Let $\A^{p}_{g,\al}$ be the $\C$-vector subspace of $\Gamma(U,\Omega_{U}^{p}\,g^{\al})$ generated by
$$g^{\al}\1\om_{k_1}\sw\cdots\sw\om_{k_{p}}\q\h{for any}\q k_1<\cdots<k_{p},$$
(see \cite[Lemma 5]{Br} for the case $\al=0$). The complex $\A^{\ssb}_{g,\al}$ with differential given by $\oma{\wedge}\,$ is identified with a subcomplex of $\Gamma(U,\Omega_{U}^{\ssb}\,g^{\al})$, that is, there is a canonical morphism of complexes
$$(\A^{\ssb}_{g,\al},\oma{\wedge})\to\Gamma(U,\Omega_{U}^{\ssb}\,g^{\al}).
\leqno(2.1.3)$$
\sk
In the case $\al=0$ (that is, $\al_k=0$ ($k\in[1,d]$) so that $g^{\al}=1$, we denote $\A^{\ssb}_{g,\al}$ by $\A^{\ssb}$. This is called the {\it Orlik-Solomon algebra,} see \cite{OS}. There is a trivial isomorphism
$$(\A^{\ssb},\oma{\wedge})=(\A^{\ssb}_{g,\al},\oma{\wedge}),$$
(since the exterior product with $\oma$ commutes with the multiplication by $g^{\al}$). We thus get the canonical morphism of complexes
$$(\A^{\ssb},\oma{\wedge})\to\Gamma(U,\Omega_{U}^{\ssb}\,g^{\al}).
\leqno(2.1.4)$$
The source is called the {\it Aomoto complex} associated with $\al\in\C^{d-1}$. By Corollary~(1.5), we have the canonical isomorphisms
$$\A^p\simto H^p(U,\C)=H^0\bl(Y,\Om_Y^p(\lH X)\br)\q(p\ges 0),
\leqno(2.1.5)$$
where the first isomorphism is due to Brieskorn \cite[Lemma 5]{Br}. Put
$$\al_d=-\msum_{k=1}^{d-1}\,\al_k,$$
so that
$$\msum_{k=1}^d\,\al_k=0.
\leqno(2.1.6)$$
\sk
We have the following.
\msn
{\bf Theorem~2.2} (\cite{ESV}, \cite{STV}). {\it The morphism $(2.1.4)$ is a quasi-isomorphism if condition~$(3)$ in the introduction is satisfied for any {\it dense edge} $E$ of $X$.}
\msn
{\bf Remarks.} (i) In the case of a constant local system with $\al=0$, this is due to Brieskorn \cite[Lemma 5]{Br} as in (2.1.5). Under a condition stronger than (2.1.7), the quasi-isomorphisms (2.1.3--4) were shown in \cite{ESV}, and it was shown in \cite{STV} that it is enough to assume condition (2.1.7) only for {\it dense edges}.
\ms
(ii) We say that an edge of a hyperplane arrangement $X$ is {\it dense} if the hyperplanes in $X$ containing the edge form an {\it indecomposable} arrangement. An arrangement is called {\it indecomposable} if its defining polynomial is not a product of two (non-constant) polynomials with disjoint variables, see \cite{STV} for more details.
\ms
(iii) A local system $L$ of rank 1 on $U$ is determined (up to a non-canonical isomorphism) by the eigenvalues of its local monodromies $\la_k$. This means that local systems $L$ and $L'$ of rank 1 are (non-canonically) isomorphic if the eigenvalues of their local monodromies coincide. Indeed, the rank 1 local system $L^{\vee}\otimes L'$ can be extended over $Y$, and it is trivial, since $Y$ is simply connected. Here $L^{\vee}$ is the dual line bundle of $L$.
\msn
{\bf 2.3.~Cohomology of the complement.} It is rather easy to calculate the cohomology of the complement $H^{\ssb}(U,\Q)$, or its dual $H_c^{\ssb}(U,\Q)$, using the weight filtration $W$ of the mixed Hodge module $(j_U)_!\Q_{h,U}[n{-}1]$ as follows. Here $j_U:U\into Y$ is the inclusion, and $\Q_{h,U}[n{-}1]$ denotes in this paper the pure Hodge module of weight $n{-}1$ with constant coefficients on $U$, see \cite{mhm}. This also clarifies the {\it geometric} meaning of the M\"obius function studied in \cite{OS}.
\sk
By the {\it strict support decomposition} of pure Hodge modules (see \cite[5.1.3.5]{mhp}), there are pure Hodge modules $\M_{Z,k}$ of weight $k$ and with strict support $Z\in\Sc_X$ in the notation of (1.1) such that
$$\Gr^W_k(j_U)_!(\Q_{h,U}[n{-}1])=\mopl_{Z\in\Sc_X}\,\M_{Z,k}.$$
Here $\M_{Z,k}$ is a {\it constant} pure Hodge module on $Z$, since it has the strict support $Z$ (that is, its underlying $\Q$-complex is an intersection complex with some local system coefficients \cite{BBD}), and is constant over a dense open subvariety of $Z$ (considering the arrangement consisting of all the hyperplanes in $X$ containing $Z$).
\sk
Let $i_x:\{x\}\into Y$ be the inclusion for $x\in X$. We have the following equalities in the {\it Grothendieck group of mixed Hodge modules}\,:
$$\msum_{k\in\Z}\msum_{Z\in\Sc_X}\,\bl[i_x^*\M_{Z,k}\br]=\bl[i_x^*(j_U)_!\Q_{h,U}[n{-}1]\br]=0.
\leqno(2.3.1)$$
These observations imply by decreasing induction on $d_Z:=\dim Z$ that there are (non-canonical) isomorphisms
$$\M_{Z,k}\cong\begin{cases}\Q_{h,Z}^{\rho_Z}[d_Z]&\h{if}\,\,\,k=d_Z,\\\,0&\h{if}\,\,\,k\ne d_Z,\end{cases}$$
where $\Q_{h,Z}^{\rho_Z}$ means the direct sum of $\rho_Z$ copies of $\Q_{h,Z}$. Moreover the rank $\rho_Z$ can be determined by decreasing induction on $d_Z$ using the relations coming from (2.3.1):
$$\msum_{Z'\supset Z}(-1)^{d_{Z'}-d_Z}\rho_{Z'}=0\q\q(Z\in\Sc_X),
\leqno(2.3.2)$$
(where $\rho_{Z'}=0$ for $Z'\notin\Sc_X$). This means that the integers
$$(-1)^{\gamma_Z}\rho_Z\q\h{for}\q Z\in\Sc_X$$
can be identified with the M\"obius function in \cite[Section 1.2]{OS}, where $\gamma_Z:={\rm codim}_YZ$. (This simplifies some arguments in \cite[Sections 1.7--9]{BS}.)
\sk
By duality, we then get the (non-canonical) isomorphisms
$$\Gr^W_{n-1+k}\R(j_U)_*(\Q_U[n{-}1])\cong\mopl_{Z\in\Sc_X^k}\,\Q^{\rho_Z}_Z(-k)[d_Z]\q(k\ges 0).
\leqno(2.3.3)$$
\sk
In the {\it affine} arrangement case (that is, the ambient space $Y$ is affine space), we have the $E_1$-degeneration of the following spectral sequence defined in the category of mixed $\Q$-Hodge structures:
$$E_1^{-k,j+k}=H^{j-n+1}\bl(Y,\Gr^W_{n-1+k}\1\R(j_U)_*(\Q_U[n{-}1])\br)\Longrightarrow H^j(U,\Q).
\leqno(2.3.4)$$
Indeed, the edges $Z$ are also {\it affine spaces} so that
$$H^{j-n+1}\bl(Y,\Gr^W_{n-1+k}\1\R(j_U)_*(\Q_U[n{-}1])\br)\cong\begin{cases}\mopl_{Z\in\Sc_X^k}\,\Q^{\rho_Z}_Z(-k)&\h{if}\,\,\,k=j,\\\,0&\h{if}\,\,\,k\ne j.\end{cases}
\leqno(2.3.5)$$
\sk
We thus get in the affine arrangement case:
$$H^j(U,\Q)=\mopl_{Z\in\Sc_X^j}\,\Q^{\rho_Z}(-j)\q\q(j\in[0,n{-}1]),
\leqno(2.3.6)$$
This implies that the $H^j(U,\C)$ are generated by integral logarithmic forms as in Remark~(i) after Corollary~(1.5), see (2.4) below for a more precise assertion.
\sk
In the projective arrangement case, we have (2.3.6) with $\Sc_X^j$ replaced by $\Sc_{X'}^j$, where $X'$ is the affine arrangement defined as in (2.1) by
$$X':=X\setminus X_d\,\,\subset\,\,Y\setminus X_d=\C^{n-1}.$$
\msn
{\bf Remarks.} (i) Assume $n=3$ in the notation of (2.1) for simplicity. The relations among the $\om_j\sw\om_k$ are generated by
$$\om_i\sw\om_j+\om_j\sw\om_k+\om_k\sw\om_i=0\q\h{for}\q X'_i\cap X'_j\cap X'_k\ne\emptyset,
\leqno(2.3.7)$$
(see \cite{OS,OT}). This can be verified by using (2.3.2), (2.3.6). Indeed, (2.3.2) implies that
$$\rho_p=\mult_pX'-1\q\h{for}\q\{p\}\in\Sc_{X'}^2,
\leqno(2.3.8)$$
where $\mult_pX'=\#\bl\{k\in[1,d{-}1]\mid X'_k\ni p\br\}$ (since $X'$ is reduced). The relations in (2.3.7) imply the inequality
$$\dim H^2(U,\C)\les\msum_{\{p\}\in\Sc_{X'}^2}\,(\mult_pX'-1),$$
since $H^2(U,\C)$ is generated by the integral logarithmic forms $\om_j\sw\om_k$ as in \cite[Lemma 5]{Br}, see also (2.4) below. (Indeed, these forms with $X'_j\cap X'_k=\{p\}$ span a $\C$-vector space of dimension at most $\mult_pX'-1$ by using (2.3.7) for each singular point $p$ of $X'$.) We would get the strict inequality $<$ if there were more relations among the $\om_j\sw\om_k$.
\sk
So we get the direct sum decomposition
$$H^2(U,\C)=\mopl_{\{p\}\in\Sc_{X'}^2}\,H^2(U,\C)_p,
\leqno(2.3.9)$$
where $H^2(U,\C)_p$ is generated by the $\om_j\sw\om_k$ with $X'_j\cap X'_k=\{p\}$. (This is closely related to \cite[Lemma 3]{Br}.) We denote the projection $H^2(U,\C)\to H^2(U,\C)_p$ by $\pi_p$.
\ms
(ii) In the notation and assumption as in Remark (i) above, assume
$$\pi_p(\om^{\beta}\sw\oma)=0\q\h{in}\q H^2(U,\C)_p,$$
with $\om^{\beta}=\msum_k\,\beta_k\om_k$ ($\beta_k\in\C$). It is well-known (and is easy to show) that
$$\al_p\1\beta_k=\beta_p\1\al_k\q\h{if}\,\,\,\,X'_k\ni p,
\leqno(2.3.10)$$
where $\beta_p:=\msum_{X'_k\ni p}\,\beta_k$, see for instance \cite[Lemma 3.1]{LY} (or \cite[Lemma 1.4]{BDS}).
\sk
In the double point case (that is, if $\mult_pX'=2$), the above assumption trivially implies that
$$\al_j\1\beta_k=\beta_j\1\al_k\q\h{for}\,\,\,\,X'_j\cap X'_k=\{p\}.
\leqno(2.3.11)$$
\msn
{\bf 2.4.~Proof of Brieskorn's results using Hodge theory.} It is rather easy to prove the assertion in \cite[Lemma 5]{Br} by induction on the dimension $n$ and the number of hyperplanes $d$ using the mixed Hodge theory as follows:
\sk
For an {\it affine} hyperplane arrangement $X=\mcup_{k=1}^dX_k$ in $Y=\C^n$, Set
$$\aligned&X':=\mcup_{k=1}^{d-1}X_k,\q X'':=X_d\cap X'\subset X_d,\\&U:=Y\setminus X,\q U':=Y\setminus X',\q U'':=X_d\setminus X'',\endaligned$$
with inclusions
$$j:U\into Y,\q j':U'\into Y,\q j'':U''\into X_d,\q i:X_d\into Y.$$
\sk
Taking the dual of the short exact sequence
$$0\to j_!\Z_U\to j'_!\1\Z_{U'}\to i_*j''_!\1\Z_{U''}\to 0,$$
and using the isomorphism
$$\D\Z_U=\Z_U(n)[2n],$$
with $\D$ the dual functor, we get the following distinguished triangle in $D^b_c(Y,\Z)$:
$$i_*\R j''_*\1\Z_{U''}(-1)[-2]\to\R j'_*\1\Z_{U'}\to\R j_*\1\Z_U\buildrel{+1}\over\to.
\leqno(2.4.1)$$
Its scalar extension by $\Z\into\Q$ can be defined in the derived category of mixed Hodge modules, and it induces a long exact sequence of mixed $\Z$-Hodge structures
$$H^{j-2}(U'')(-1)\buildrel{\al_j\,\,}\over\to H^j(U')\buildrel{\beta_j\,\,}\over\to H^j(U)\buildrel{\gamma_j\,\,}\over\to H^{j-1}(U'')(-1)\buildrel{\al_{j+1}}\over\longrightarrow H^{j+1}(U'),
\leqno(2.4.2)$$
where the cohomology is with $\Z$-coefficients. By induction on $n$ and $d$, we see that the $H^j(U)$ are torsion-free, and have pure weight $2j$ and type $(j,j)$, hence the $\al_j$ vanish, the $\beta_j$ are injective, and the $\gamma_j$ are surjective (all with $\Z$-coefficients).
\sk
We can moreover show that the cohomology groups $H^j(U)$ are generated by the classes of exterior products of the $\ddd g_k/g_k$ (up to a Tate twist, see \cite{De}) with $g_k$ a defining polynomial of $X_k$ with degree 1, and the morphism $\gamma_j$ for these forms is obtained by taking the residue along $X_d$, that is, by dividing out the forms by $\ddd g_d/g_d$ (if divisible, and it vanishes otherwise). Indeed, these assertions can be reduced to the normal crossing case by considering hyperplane arrangements with normal crossings $X^{\rm nc}$ contained in $X$, and using the injectivity of $\beta_j$ for the inclusion $U\into Y\setminus X^{\rm nc}$ by factorizing the latter. In the case $X$ is a divisor with normal crossings on $Y$, we have the short exact sequences as is well-known:
$$0\to\Om_Y^p(\log X')\to\Om_Y^p(\log X)\to\Om_{X_d}^{p-1}(\log X'')\to 0\q(p\in\Z),
\leqno(2.4.3)$$
where the last morphism is given by the residue along $X_d$. This short exact sequence can be obtained also by applying the filtered de Rham functor (see \cite{mhp}) and $\Gr_F^p$ to the short exact sequence of filtered regular holonomic $\DD_Y$-modules corresponding to (2.4.1):
$$0\to j'_*(\Oc_{U'},F)\to j_*(\Oc_U,F)\to i_*j''_*(\Oc_{U''},F[-1])\to0.
\leqno(2.4.4)$$
\sk
We can also prove \cite[Lemma 3]{Br} by using the above assertions together with (2.3.6). (See \cite{Br} and also \cite{JR}, \cite{OT} for arguments without using the mixed Hodge theory.)
\msn
{\bf 2.5.~Calculation of the Milnor cohomology.} Let $f$ be a defining polynomial of a projective hyperplane arrangement in $\PP^{n-1}$. Set $F_{\!f}:=f^{-1}(1)\subset\C^n$. This is the Milnor fiber of $f$. It is well-known (see for instance \cite{Di2}, \cite[Sections 1.3--4]{BS}) that the $\la$-eigenspace of the monodromy on the Milnor cohomology $H^j(F_{\!f},\C)_{\la}$ is calculated by the cohomology of a rank 1 local system such that the local monodromies around any hyperplanes are $\la^{-1}$, and moreover $H^j(F_{\!f},\C)_{\la}=0$ unless $\la^d=1$ with $d:=\deg f$ (using the geometric monodromy for a homogeneous polynomial case). As a corollary, we get (see for instance \cite[1.4.2]{BS})
$$\msum_j\,(-1)^j\dim H^j(F_{\!f},\C)_{\la}=\chi(U)\q\h{if}\q\la^d=1.
\leqno(2.5.1)$$
\sk
Set $\la=\exp(-2\pi ik/d)$ for $k\in[1,d-1]$ following \cite[Sections 3.2--3]{BSY} instead of \cite{BDS}. (Note that $\la\ne 1$.) We can calculate the monodromy eigenspaces of the Milnor cohomology groups $H^j(F_{\!f},\C)_{\la}$ by using the Aomoto complex as in \cite{ESV,STV} (see Theorem~(2.2) above) if there is a subset $I\subset\{1,\dots,d\}$ with $|I|=k$ and such that condition~(3) in the introduction is satisfied by setting
$$\al_j:=\begin{cases}1-\kod&\h{if}\q j\in I,\\-\kod&\h{if}\q j\notin I.\end{cases}
\leqno(2.5.2)$$
(Note that the sign should be reversed if we set $\la=\exp(2\pi ik/d)$ as in \cite{BDS}.)
\msn
{\bf Examples 2.6.} (i) There are many examples such that the $H^j(F_{\!f},\C)_{\la}$ cannot be calculated by the method in (2.5) when $d\ges 9$. For instance, let
$$f=xy(x-2y)(y+z)(2x-4y+z)(x-y+z)(2x-y+z)(x+y+z)z.$$
The picture of the affine part of the arrangement is obtained by setting $z=1$ as follows:
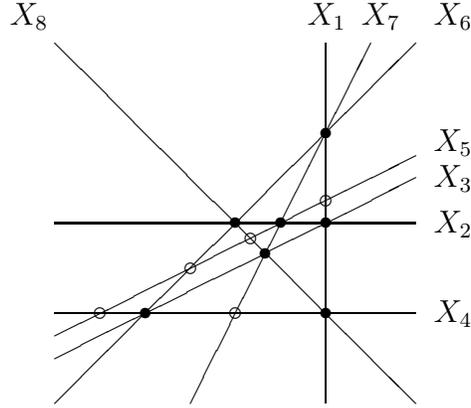
\sk
$$\setlength{\unitlength}{1.2cm}
\begin{picture}(4,4.3)
\put(-0.5,4.2){$X_8$}
\put(2.8,4.2){$X_1$}
\put(3.4,4.2){$X_7$}
\put(4.2,4.2){$X_6$}
\put(4.2,2.8){$X_5$}
\put(4.2,2.4){$X_3$}
\put(4.2,1.9){$X_2$}
\put(4.2,0.9){$X_4$}
\put(3,0){\line(0,1){4}}
\put(0,2){\line(1,0){4}}
\put(0,1){\line(1,0){4}}
\put(0,0.5){\line(2,1){4}}
\put(0,4){\line(1,-1){4}}
\put(0,0){\line(1,1){4}}
\put(1.5,0){\line(1,2){2}}
\put(0,0.75){\line(2,1){4}}
\put(3,1){\circle*{.12}}
\put(3,2){\circle*{.12}}
\put(3,3){\circle*{.12}}
\put(2.5,2){\circle*{.12}}
\put(2,2){\circle*{.12}}
\put(2.33,1.67){\circle*{.12}}
\put(1,1){\circle*{.12}}
\put(0.5,1){\circle{.12}}
\put(1.5,1.5){\circle{.12}}
\put(2.17,1.83){\circle{.12}}
\put(2,1){\circle{.12}}
\put(3,2.25){\circle{.12}}
\end{picture}$$
\skn
The hyperplane defined by the $k$\1th factor of $f$ is denoted by $X_k$ ($k\in[1,8]$) except for the line at infinity $X_9=\{z=0\}$.
We denote double and triple points by white and black vertices respectively. Note that parallel lines meet at infinity giving a triple point of the projective arrangement $X$.
\sk
We see that there is no subset $I\subset\{1,\dots,9\}$ with $|I|=3$ and such that $X_k\cap X_j$ is a double point of the projective arrangement $X$ for any $k,j\in I$ with $k\ne j$. This means that condition~(3) cannot be satisfied by defining the $\al_j$ as in (2.5.2).
\sk
The Milnor cohomology $H^1(F_{\!f},\C)_{\la}$ for $\la=\exp(-2\pi i/3)$ cannot be calculated by applying Theorem~1 either, since we need a modification as follows: In this case the subset $I$ of $\{1,\dots,9\}$ with $|I|=3$ must be $\{4,5,9\}$ (since $X_4\cap X_5$ must be a double point of $X$), and
$$I^b=\{X_2,X_3,X_4,X_5\},\q I^c=\{X_1,X_6,X_7,X_8,X_9\}.$$
We have the connectivity of
$$\aligned&(X_1\cup X_6\cup X_7)\setminus(X^{I^c(0)}\cup X^{I^b}),\\&(X_1\cup X_5\cup X_6\cup X_7\cup X_8)\setminus(X^{I^c(0)}\cup X^{I^b}),\endaligned$$
and $X_5$ is good for the first one, that is,
$$\#\bl(X_5\cap(X_1\cup X_6\cup X_7)\setminus(\mcup_{j\in I^b_{<5}}\,X_j)\br)\ges 2.$$
\ms
(ii) As a more complicated example with $d=9$, we have
$$f=(x^3-y^3)(y^3-z^3)(x^3-z^3).$$
There are {\it no double points} on $X$, although it has $12$ {\it triple points}. We cannot apply the method in this paper because of condition~(d) in Theorem~1.

\bs\bs
\vbox{\centerline{\bf 3. Three variable case}
\bsn
In this section we restrict to the case $n=3$, and prove Theorem~1 for the calculation of the cohomology of a rank 1 local system via the Aomoto complex.}
\msn
{\bf Proposition~3.1.} {\it Assume $n=3$, and $d:=\deg X\ges 3$. Let $C\subset Y=\PP^2$ be a general line. In the notation of $(1.3)$, we have for any $m\ges0:$
$$H^0\bl(C,\Om_Y^p(\lH X)(m)|_C\br)=\begin{cases}m{+}1&\h{if}\,\,\,\,p=0,\\2m{+}d{-}1&\h{if}\,\,\,\,p=1,\\m{+}d{-}2&\h{if}\,\,\,\,p=2,\end{cases}
\leqno(3.1.1)$$
$$H^j(C,\Om_Y^p(\lH X)(m)|_C\br)=0\q\q(j>0,\,p\in\N),
\leqno(3.1.2)$$
where $(m)$ means the tensor product with $\Oc_Y(m)$, and $|_C$ is the restriction as $\Oc$-modules, that is, the tensor product with $\Oc_C$ over $\Oc_Y$.}
\msn
{\it Proof.} For $p=0$, the assertions follows from Remark (ii) after Corollary~(1.5).
\sk
For $p=1$, we first show the short exact sequence
$$0\to\Nc^*_{C/Y}(m)\to\Om^1_Y(\lH X)(m)|_C\to\Om_C^1(X_C)(m)\to 0.
\leqno(3.1.2)$$
Here $\Nc^*_{C/Y}$ is the conormal sheaf of $C\subset Y$, $(X_C)$ means the tensor product of $\Oc_C(X_C)$ over $\Oc_C$ with $X_C:=X\cap C$, and $X$ is viewed as a {\it reduced\1} divisor on $Y$. Since $\Oc_Y(m)|_C=\Oc_C(m)$ (with $Y=\PP^2$, $C=\PP^1$), we may assume $m=0$. There is a canonical morphism
$${\rm Ker}(\Om^1_Y|_C\to\Om_C^1)\to{\rm Ker}\bl(\Om^1_Y(\lH X)|_C\to\Om_C^1(X_C)\br).
\leqno(3.1.4)$$
We see that this is an isomorphism using local coordinates, since $C$ intersects the projective line arrangement $X$ transversally at smooth points. The exact sequence (3.1.2) then follows.
\sk
We now get the assertions for $p=1$ by (3.1.2), since $C=\PP^1$, $d\ges 3$, $m\ges 0$, and
$$\Nc^*_{C/Y}=\Oc_C(-1),\q\Om_C^1(X_C)=\Oc_C(d-2).
\leqno(3.1.5)$$
(Note that $\dim H^j(C,\Oc_C(k))=k+1$ for $j=0$, and $0$ otherwise, assuming $k\ges -1$.)
\sk
As for the case $p=2$, we have the isomorphism
$$\Om^2_Y(\lH X)(m)|_C=\Om_C^1(X_C)\otimes_{\Oc_C}\Nc^*_{C/Y}(m),
\leqno(3.1.6)$$
since $\Om^2_Y(\lH X)$ is isomorphic to
$$\Om^2_Y(X):=\Om_Y^2\otimes_{\Oc_Y}\Oc_Y(X)$$
on a neighborhood of $C\subset Y$. So the assertions for $p=2$ also follow from (3.1.5) since $d\ges 3$. This finishes the proof of Proposition~(3.1).
\msn
{\bf Corollary~3.2.} {\it In the assumptions of Proposition~$(3.1)$, we have the following surjective canonical morphisms for any $\,m,i\ges 0$, $p\in\Z:$
$$R_{Y,i}\otimes_{\C}H^0\bl(Y,\Om^p_Y(\lH X)(m)\br)\to H^0\bl(Y,\Om^p_Y(\lH X)(m{+}i)\br),
\leqno(3.2.1)$$
where $R_{Y,i}:=H^0\bl(Y,\Oc_Y(i)\br)$.}
\msn
{\it Proof.} We may assume $i=1$. Let $x,y,z$ be a $\C$-basis of $R_{Y,1}$ such that $z^{-1}(0)=C$ in the notation of Proposition~(3.1). The morphism
$$z\otimes:H^0\bl(Y,\Om^p_Y(\lH X)(m)\br)\into H^0\bl(Y,\Om^p_Y(\lH X)(m{+}1)\br)$$
is identified with the inclusion
$$H^0\bl(Y,\Om^p_Y(\lH X)(-C)(m{+}1)\br)\into H^0\bl(Y,\Om^p_Y(\lH X)(m{+}1)\br),$$
where $(-C)$ means the tensor product over $\Oc_Y$ with $\Oc_Y(-C)$, that is, the ideal sheaf of $C$.
\sk
Let $x',y'$ be the restrictions of $x,y$ to $R_{C,1}:=H^0\bl(C,\Oc_C(1)\br)$. Since these are generators, the surjectivity of (3.2.1) for $i=1$ is then reduced to the following surjections for $p,m\ges 0:$
$$R_{C,1}\otimes_{\C}H^0\bl(C,\Om^p_Y(\lH X)|_C(m)\br)\onto H^0\bl(C,\Om^p_C(\lH X)|_C(m{+}1)\br),
\leqno(3.2.2)$$
using the isomorphisms between the cokernels of the above two inclusions, if we have the following canonical surjections for $p,m\ges 0:$
$$H^0\bl(Y,\Om^p_Y(\lH X)(m)\br)\onto H^0\bl(C,\Om^p_Y(\lH X)|_C(m)\br).
\leqno(3.2.3)$$
Here $m$ must be $\ges 0$, and not $\ges 1$, since we use the compatibility of these surjections with the actions of $x,y$ and $x',y'$.
\sk
The surjectivity of (3.2.3) follows from Corollary~(1.4) if $m\ges 1$. This holds also for $m=0$, since the vector spaces
$$H^0\bl(C,\Om^p_Y(\lH X)|_C\br)\q\h{for}\,\,\,\,p=1\,\,\,\h{and}\,\,\,2$$
are generated respectively by the images of 
$$\om_k\,\,(k\in[1,d{-}1])\q\h{and}\q\om_1\sw\om_k\,\,(k\in[2,d{-}1]),$$
where the $\om_k$ are as in (2.1), and are logarithmic also along $X_d\subset Y$. Indeed, the above images are linearly independent in the target vector spaces by taking the residues along $C\cap X_j$ ($j\in[1,d{-}1]$) after restricting them to $C$ as differential forms if $p=1$ (see also (3.1.2)), and looking at their poles along $C\cap X_j$ after restricting them to $C$ as meromorphic sections of a locally free sheaf on $Y$ if $p=2$. Note that the dimensions of the target vector spaces are given in Proposition~(3.1). (The case $p=0$ is trivial.)
\sk
The surjectivity of (3.2.2) can be shown by using (3.1.2), (3.1.5--6) (together with the snake lemma) since $d\ges 3$. This finishes the proof of Corollary~(3.2).
\msn
{\bf Proposition~3.3.} {\it In the notation of $(2.1)$ and the assumption of Proposition~$(3.1)$, set
$$I^{(k)}:=\{j\in[1,d{-}1]\mid X_j\supset X_k\cap X_d\}\q(k\in[1,d{-}1]),$$
and assume $|I^{(1)}|=1$. Define the pole order filtration $P$ along $X_d$ by
$$P_mH^0\bl(Y,\Om^p_Y(\lH X)(*X_d)\br):=H^0\bl(Y,\Om^p_Y(\lH X)(mX_d)\br)\q(m\ges 0).$$
Here $P_m=0$ for $m<0$, and $(mX_d)$ means the tensor product with $\Oc_Y(mX_d)$ over $\Oc_Y$. Set $x:=g_1$, $y:=g_2$ so that $\om_1=\tfrac{\dx}{x}$, $\om_2=\tfrac{\dy}{y}$. Then, for $m\ges 1$, the vector spaces
$$\Gr^P_mH^0\bl(Y,\Om^p_Y(\lH X)(*X_d)\br)$$
have $\C$-bases consisting of the classes of}
$$x^m\1\om_k\,\,\,(k\in[3,d{-}1]),\,\,\,x^iy^{m-i}\1\tfrac{\dx}{x},\,\,\,x^iy^{m-i}\1\tfrac{\dy}{y}\,\,\,(i\in[0,m])\q\h{if}\q p=1,
\leqno(3.3.1)$$
$$x^{m-1}\dx{\sw}\om_k\,\,\,(k\in[3,d{-}1]),\,\,\,x^iy^{m-i}\1\tfrac{\dx}{x}{\sw}\tfrac{\dy}{y}\,\,\,(i\in[0,m])\q\h{if}\q p=2.
\leqno(3.3.2)$$
\msn
{\it Proof.} For $p=1,2$, set
$$\Psi^p:=H^0\bl(Y,\Om^p_Y(\lH X)(*X_d)\br).$$
Using the isomorphisms $\Oc_Y(mX_d)\cong\Oc_Y(m)$ together with Corollary~(1.4), we get the short exact sequences for $m\ges 1$:
$$0\to P_{m-1}\Psi^p\to P_m\Psi^p\to H^0\bl(C,\Om^p_Y(\lH X)(m)|_{C}\br)\to 0.$$
Combined with Proposition~(3.1), this implies that
$$\dim\Gr^P_m\Psi^p=\begin{cases}2m+d-1&\h{if}\,\,\,\,p=1,\\m+d-2&\h{if}\,\,\,\,p=2.\end{cases}$$
So it is enough to show that the $\Gr^P_m\Psi^p$ are generated by the classes in (3.3.1--2) for $p=1,2$.
\sk
By Corollary~(3.2), we have the surjections
$$\Gr^P_mH^0\bl(Y,\Oc_Y(*X_d)\br)\otimes_{\C}H^0\bl(Y,\Om^p_Y(\lH X)\br)\onto\Gr^P_m\Psi^p\q(m\ges 1),
\leqno(3.3.3)$$
since $\Oc_Y(mX_d)\cong\Oc_Y(m)$. Note that the $g_k$ have poles of order 1 along $X_d$, and
$$\Gr^P_mH^0\bl(Y,\Oc_Y(*X_d)\br)=\msum_{i=0}^m\,\C\1[x^iy^{m-i}],
\leqno(3.3.4)$$
since $X'_2$ is not parallel to $X'_1$ in $\C^2$ (that is, $X'_1\cap X'_2\ne\emptyset$) by the assumption: $|I^{(1)}|=1$. Here the pole order filtration $P$ on $H^0\bl(Y,\Oc_Y(*X_d)\br)$ is defined similarly to the one on $\Psi^p$, and $[v]$ means the class of $v$ in $\Gr^P_m$ in general.
\sk
On the other hand we have by (2.1.5)
$$H^0\bl(Y,\Om^p_Y(\lH X)\br)=\begin{cases}\sum_{k=1}^{d-1}\,\C\1\om_k&\h{if}\,\,\,p=1,\\ \sum_{1\les j<k\les d-1}\,\C\1\om_j\sw\om_k&\h{if}\,\,\,p=2.\end{cases}
\leqno(3.3.5)$$
\sk
For $k,i,j\in[1,d{-}1]$ with $\om_i\sw\om_j\ne 0$, or equivalently, $X'_i\cap X'_j\ne\emptyset$ (that is, $X'_i$ is not parallel to $X'_j$), there are complex numbers $c_{k,i,j},c'_{k,i,j},c''_{k,i,j}$ satisfying
$$g_k=c_{k,i,j}+c'_{k,i,j}\,g_i+c''_{k,i,j}\,g_j\,\,\,\,(\h{hence}\,\,\,\ddd g_k=c'_{k,i,j}\,\ddd g_i+c''_{k,i,j}\,\ddd g_j).
\leqno(3.3.6)$$
Setting $i=1$, $j=2$ in (3.3.6), we get the following equalities in $P_m\Psi^1$ for $k\in[3,d{-}1]$ and $i,i'\in\N$ with $i+i'=m-1$\,:
$$\aligned&x^iy^{i'}(c_{k,1,2}+c'_{k,1,2}\1x+c''_{k,1,2}\1y)\om_k=x^iy^{i'}\dg_k\\&=x^iy^{i'}(c'_{k,1,2}\1\dx+c''_{k,1,2}\1\dy).\endaligned
\leqno(3.3.7)$$
Here $c''_{k,1,2}\ne 0$, since $|I^{(1)}|=1$ (and hence $X'_k$ is not parallel to $X'_1$). These imply by decreasing induction on $i\in[0,m{-}1]$ that $[x^iy^{m-i}\om_k]$ belongs to the subspace of $\Gr^P_m\Psi^1$ generated by the terms in (3.3.1) if $k\in[3,d{-}1]$. So the assertion for $p=1$ follows.
\sk
For $p=2$, we have the following equalities in $P_1\Psi^2$ by (3.3.6) if $\om_i\sw\om_j\ne 0\,$:
$$\aligned g_k\,\om_i\sw\om_j&=c_{k,i,j}\,\om_i\sw\om_j+c'_{k,i,j}\,\dg_i\sw\om_j-c''_{k,i,j}\,\dg_j\sw\om_i\\
&=c_{k,i,j}\,\om_i\sw\om_j+\dg_k\sw\om_j-\dg_k\sw\om_i\q(k=1,2),\endaligned
\leqno(3.3.8)$$
$$(c'_{j,1,2}\1\dx+c''_{j,1,2}\1\dy){\sw}\om_j=\dg_j{\sw}\om_j=0\,\,\,\,\h{(similarly for $\om_i$)},
\leqno(3.3.9)$$
with $c''_{j,1,2}\ne 0$ for $j\ne 1$ (since $|I^{(1)}|=1$). These imply that $\Gr^P_m\Psi^2$ is generated by the classes of
$$x^iy^{i'}\tfrac{\dx}{x}\sw\tfrac{\dy}{y},\q x^iy^{i''}\dx\sw\om_k\,\,\,(k\in[3,d{-}1]),$$
where $i,i',i''$ are non-negative integers with $i+i'=m$, $i+i''=m-1$.
\sk
By (3.3.7) (with $m$ replaced by $m-1$), we get the following equalities in $P_m\Psi^2$:
$$x^iy^{i'}(c_{k,1,2}+c'_{k,1,2}\1x+c''_{k,1,2}\1y)\dx\sw\om_k=c''_{k,1,2}\,x^iy^{i'}\1\dx\sw\dy.
\leqno(3.3.10)$$
Here $i,i'$ are non-negative integers with $i+i'=m-2$. The assertion for $p=2$ then follows by an inductive argument similar to the case $p=1$. This finishes the proof of Proposition~(3.3).
\msn
{\bf 3.4.~Proof of Theorem~1.} In the notation of the proof of Proposition~(3.3), we have a short exact sequence of complexes
$$0\to P_0\Psi^{\ssb}\to P_1\Psi^{\ssb}\to\Gr^P_1\Psi^{\ssb}\to 0,$$
where the twisted differential $\ddd^{\al}$ of the filtered complex $\Psi^{\ssb}$ is defined by using $g^{\al}$ as in (2.1). It induces the long exact sequence
$$\to H^jP_0\Psi^{\ssb}\to H^jP_1\Psi^{\ssb}\to H^j\Gr^P_1\Psi^{\ssb}\to H^{j+1}P_0\Psi^{\ssb}\to.$$
Here $P_0\Psi^{\ssb}$ is the Aomoto complex by (2.1.5), and we have the isomorphisms
$$H^jP_1\Psi^{\ssb}\simto H^j(L,\C)\q(j\in\Z),$$
by condition~(a) together with Proposition~(1.2) and Corollary~(1.4), generalizing some arguments in \cite{ESV} as in \cite{bha}. To show Theorem~1, it is then enough to show the injectivity of
$$H^1\Gr^P_1\Psi^{\ssb}\to H^2P_0\Psi^{\ssb}.$$
However, we will use another connecting morphism which is easier to treat, see (3.4.1) below.
\sk
By condition~(d), we may assume
$$|I^{(1)}|=|I^{(2)}|=1,$$
by changing the order of $\{1,\dots,d\}$ if necessary. (So the assumption of Proposition~(3.3) is satisfied.) We then see that $\Gr^P_1\Psi^1$, $\Gr^P_1\Psi^2$ are respectively generated over $\C$ by
$$x\1\om_k\,\,\,(k\in[3,d{-}1]),\,\,\,\dx,\,\,\,\dy,\,\,\,x\1\oma,\,\,\,y\1\oma,$$
$$\dx{\sw}\om_k\,\,\,(k\in[3,d{-}1]),\,\,\,\dx{\sw}\oma,\,\,\,\dy{\sw}\oma,$$
and the complex $\Gr^P_1\Psi^{\ssb}$ has the acyclic subcomplex $\Theta^{\ssb}$ generated over $\C$ by the classes of
$$x,\,y,\q\dx,\,\dy,\,x\1\oma,\,y\1\oma,\q\dx\sw\oma,\,\dy\sw\oma.$$
So $H^1\Gr^P_1\Psi^{\ssb}$ can be calculated by using its quotient complex
$$\Phi''{}^{\ssb}:=\Gr^P_1\Psi^{\ssb}/\Theta^{\ssb}.$$
\sk
Since the above acyclic subcomplex $\Theta^{\ssb}$ is naturally lifted to the acyclic subcomplex $\Theta'{}^{\ssb}$ of $P_1\Psi^{\ssb}$, we can consider its quotient complex
$$\Phi^{\ssb}:=P_1\Psi^{\ssb}/\Theta'{}^{\ssb},$$ 
so that we get the short exact sequence of complexes
$$0\to\Phi'{}^{\ssb}\to\Phi^{\ssb}\to\Phi''{}^{\ssb}\to0,$$
with $\Phi'{}^{\ssb}=P_0\Psi^{\ssb}$. It is then sufficient to prove the injectivity of
$$\dd:H^1\Phi''{}^{\ssb}\to H^2\Phi'{}^{\ssb},
\leqno(3.4.1)$$
using the canonical morphism between the above two short exact sequences of complexes.
\sk
By the preceding arguments, $\Phi''{}^1$ is generated by the classes of $x\1\om_k$ ($k\in[3,d{-}1])$, and the images of these $x\1\om_k$ under the twisted differential $\ddd^{\al}$ are given by
$$\aligned\da(x\om_k)&=\dx\sw\om_k+\bl(\msum_{j=1}^{d-1}\al_j\om_j\br)\sw x\om_k\\
&=\msum_{j\in I_{p(k)}^c}\,c_{1,j,k}\al_j\om_j\sw\om_k+\bl(1+\msum_{j\in I_{p(k)}^c}\,\al_j\br)\dx\sw\om_k\\&\q-\dx\sw\bl(\msum_{j\in I_{p(k)}^c}\,\al_j\om_j\br),\endaligned
\leqno(3.4.2)$$
where $\{p(k)\}=X_k\cap X_d$, and
$$I_{p(k)}:=\{j\in[1,d]\mid X_j\ni p(k)\},\q I_{p(k)}^c:=\{1,\dots,d\}\setminus I_{p(k)}.$$
We then get
$$\aligned\da(x\om_k+x\1\oma)&=\msum_{j\in I_{p(k)}^c}\,c_{1,j,k}\al_j\om_j\sw\om_k+(1-\al_{p(k)})\dx\sw\om_k\\
&\q+\dx\sw\bl(\msum_{j\in I'_{p(k)}}\,\al_j\om_j\br)-\dx\sw\oma,\endaligned
\leqno(3.4.3)$$
using condition~(1) in the introduction (where $I'_{p(k)}:=I_{p(k)}\setminus\{d\}$). This implies for $p\in X_d$
$$\aligned&\da\bl(\msum_{k\in I'_p}\,x\al_k\om_k+\al'_px\1\oma\br)\\&\equiv\msum_{k\in I'_p}\msum_{j\in I_p^c}\,c_{1,j,l}\al_j\om_j\sw\al_k\om_k+(1-\al_d)\dx\sw\bl(\msum_{k\in I'_p}\,\al_k\om_k)\\
&\q\,\,\,({\rm mod}\,\,\,\C\,\dx\sw\oma),\endaligned
\leqno(3.4.4)$$
since $\al'_p:=\msum_{k\in I'_p}\,\al_k=\al_p-\al_d$, where $I_p$, $\al_p$ are as in the introduction, and $I'_p:=I_p\setminus\{d\}$. Note that $\al_d\ne 1$ by condition~(a).
\sk
The above calculations imply the surjectivity (and hence the injectivity) of the following subcomplex of $\Phi''{}^{\ssb}$:
$$\mopl_{k\in I'_p}\,\C[x\,\om_k]\to\mopl_{k\in I'_p}\,\C[\dx\sw\om_k]\q\h{if}\q\al_p\ne 1.
\leqno(3.4.5)$$
So any element of $H^1\Phi''{}^{\ssb}$ is represented in $\Phi^1$ by an element $\eta$ which is written as
$$\eta=\msum_{k\in I'{}^b}\,a_k\,x\1(\om_k+x\1\oma)\q\h{with}\q a_k\in\C,$$
where $I'{}^b:=I^b\setminus\{d\}$ with $I^b$ as in Theorem~1. By (3.4.3) the vanishing of its image in $\Phi''{}^2$ is determined by looking at the coefficients of $\dx\sw\bl(\msum_{j\in I'_{p_i}}\,\al_j\om_j\br)$ for each $p_i$, since $\al_{p(k)}=1$ for $k\in I^b$ by condition~(a). (Note that $p(k)=p_i$ for $k\in I_{p_i}$.) So we get the equality
$$\ddd^{\al}\eta=\msum_{k\in I'{}^b}\,a_k\bl(\msum_{j\in I^c_{p(k)}}\,c_{1,j,k}\al_j\om_j\sw\om_k\br)\q\h{in}\q\Phi'{}^2=P_0\Psi^2,
\leqno(3.4.6)$$
and the image of $[\eta]\in H^1\Phi''{}^{\ssb}$ by the connecting morphism $\dd$ of the long exact sequence associated with the short exact sequence is given by the right-hand side of (3.4.6).
\sk
Assume it vanishes in $H^2\Phi'{}^{\ssb}$. This means that we have the equality
$$\ddd^{\al}\eta=\msum_{k=1}^d\,b_k\bl(\msum_{j=1}^{d-1}\,\al_j\om_j\sw\om_k\br)\q\h{with}\q b_k\in\C.$$
Using condition (b) together with Remark~(ii) after (2.3), we see that there is $c\in\C$ satisfying the condition
$$b_k=c\al_k\q\h{for any}\,\,\,k\in I^c.
\leqno(3.4.7)$$
By the relation $\oma\sw\oma=0$ together with the decomposition
$$\oma=\msum_{j=1}^{d-1}\,\al_j\om_j=\msum_{k\in I'{}^b}\,\al_k\om_k+\msum_{k\in I^c}\,\al_k\om_k,$$
condition (3.4.7) implies the equality
$$\ddd^{\al}\eta=\msum_{k\in I'{}^b}\,b'_k\bl(\msum_{j=1}^{d-1}\,\al_j\om_j\sw\om_k\br)\q\h{with}\q b'_k\in\C.
\leqno(3.4.8)$$
Comparing (3.4.6) and (3.4.8), we then get $a_k=b'_k=0$ using condition~(c) and Remark~(i) after (2.3). (Note that $c_{1,j,k}=g_1(p_{j,k})$ with $\{p_{j,k}\}=X_j\cap X_k$, and the function $g_1$ is not constant on any $X'_k$ ($k\in I'{}^b$) since $|I^{(1)}|=1$.) This finishes the proof of Theorem~1.
\msn
{\bf Remark.} Condition~(c) in Theorem~1 can be replaced by a slightly weaker one:
\msn
{\rm (c)$'$} $\,\,\#\bl(X_k\cap X^{I^c}\setminus\mcup_{j\in I^b_{<k}}X_j\br)\ges 2\,\,\,\,(\forall\,k\in I^b)$,\q$\al_j\ne 0\,\,\,\,(\forall\,j\in I^c)$.
\msn
where $I^b_{<k}:=\{j\in I^b\mid j<k\}$. (However, this is still insufficient for Example~(2.6)(i).)
\bs\bs
\vbox{\centerline{\bf Appendix}
\bsn
In this Appendix we study some combinatorics of the intersection poset of the reflection hyperplane arrangement of type $G_{31}$.}
\msn
{\bf A.1.~Construction.}
Let $\La$ be the set of hyperplanes in the reduced projective hyperplane arrangement $X$ of type $G_{31}$ in $\PP^3$. Set
$$\Psi:=\mmu_4\sqcup\{0\}=\{\10,\,\,\pm\1 1,\,\,\pm\1 i\1\}\subset\C.$$
According to \cite{HR} (see also \cite{BDY}, \cite{DS}, \cite{MPP}), each hyperplane in $\La$ is defined in $\PP^3$ by a linear form
$$\ell_{\aaa}=\msum_{k=1}^4\,a_kx_k,$$
which is identified with $\aaa=(a_k)\in\Psi^4$. Here $\aaa$ is not uniquely determined by the hyperplane (which will be denoted by $X_{\aaa}$), and it has an ambiguity by the diagonal multiplicative action of $\mmu_4$ on $\Psi^4$. (Diagonal action means that $\xi\,\aaa=(\xi\1 a_k)$ for $\aaa=(a_k)\in\Psi^4$, $\xi\in\mmu_4$.) We have the inclusion
$$\La\into\Phi:=\Psi^4\!\1/\mmu_4.$$
For $\aaa=(a_k)\in\C^4$, set
$$\aligned \prd(\aaa)&:=\mprod_{k=1}^4\,a_k\,,\\ \Supp\,\aaa&:=\bl\{k\mid a_k\ne 0\br\}\subset\{1,\dots,4\}.\endaligned
\leqno{(\rm A.1.1)}$$
By {\it loc.\,cit.}, we have the decomposition
$$\aligned&\q\q\q\La=\La_1\sqcup\La_2\sqcup\La_3\q\q\h{with}\\&\La_1:=\bl\{[\aaa]\in\Phi\,\big|\,|\Supp\,\aaa|=1\br\},\\&\La_2:=\bl\{[\aaa]\in\Phi\,\big|\,|\Supp\,\aaa|=2\br\},\\&\La_3:=\bl\{[\aaa]\in\Phi\,\big|\,|\Supp\,\aaa|=4,\,\,\prd(\aaa)=\pm\1 1\br\},\endaligned
\leqno{(\rm A.1.2)}$$
where $[\aaa]\in\Phi$ denotes the class of $\aaa\in\Psi^4$.
The conditions in the definitions of $\La_1$, $\La_2$, $\La_3$ are independent of a representative $\aaa$ of $[\aaa]\in\Phi$. By (A.1.2) it is easy to see that
$$|\La_1|=4,\q|\La_2|=\tbinom{4}{2}\cdot 4=24,\q|\La_3|=4^3/2=32,\q|\La|=60.
\leqno{(\rm A.1.3)}$$
\sk
Let $\Gamma$ be the semi-direct product of the abelian group $\La_3$ (via the multiplication) with the symmetric group $\Sf_4$, that is,
$$\Gamma:=\La_3\rtimes\Sf_4.
\leqno{(\rm A.1.4)}$$
More precisely, this is defined by the subgroup of ${\rm GL}(4,\C)/\mmu_4$ consisting of the classes of $\gamma=(\gamma_{j,k})\in{\rm GL}(4,\C)$ with $\gamma_{j,k}\in\mmu_4$ ($\forall\,j,k$) and such that $\gamma_{j,k}=0$ for $k\ne\sigma_{\gamma}(j)$ and $\mprod_{j=1}^4\,\gamma_{j,\sigma_{\gamma}(j)}=\pm 1$, where $\sigma_{\gamma}\in\Sf_4$ depends on $\gamma$. There is a natural action of $\Gamma$ on $\La$ in a compatible way with the above decomposition of $\La$, and the action of $\Gamma$ is transitive on each $\La_j$ ($j=1,2,3$).
\msn
{\bf A.2.~Notation.} (i) For $[\aaa],[\bb]\in\La$ with $[\aaa]\ne[\bb]$, set
$$\aligned
\La(\aaa,\bb)&:=\bl\{[\cc]\in\La\,\,\big|\,\,X_{\cc}\supset X_{\aaa}\cap X_{\bb}\br\}\setminus\{[\aaa],[\bb]\},\\
\La_j(\aaa,\bb)&:=\La_j\cap\La(\aaa,\bb),\\
\lah(\aaa,\bb)&:=\La(\aaa,\bb)\sqcup\{[\aaa],[\bb]\},\\
\mult(\aaa,\bb)&:=\bl|\lah(\aaa,\bb)\br|\in\Z_{\ges 2},\\
\endaligned
\leqno{(\rm A.2.1)}$$
where $X_{\aaa}$ is the hyperplane defined by the linear form $\ell_{\aaa}$ corresponding to $\aaa$. It is known that $\,\mult(\aaa,\bb)\in\{2,3,6\}$, see for instance \cite{BDY}.
\ms
(ii) For $k\in\{1,\dots,4\}$, we denote the projection deleting the $k\1$th factor by
$$\pi^{(k)}:\Phi=\Psi^4/\mmu_4\to\Phi':=\Psi^3/\mmu_4.
\leqno{(\rm A.2.2)}$$
\ms
(iii) For $[\aaa],[\bb]\in\La_3$, set
$$\aligned{\rm diff}(\aaa,\bb)&:=|\Supp(\aaa-\bb)|,\\{\rm diff}([\aaa],[\bb])&:=\min_{\xi\in\mmu_4}{\rm diff}(\xi\,\aaa,\bb)\in[0,3].\endaligned
\leqno{(\rm A.2.3)}$$
\ms
(iv) For $\,\aaa,\bb\in\C^4\,$ with $\,|\Supp\,\aaa|=2$, set
$$\dett(\aaa,\bb):=\dett\begin{pmatrix}a_j&a_k\\ b_j&b_k\end{pmatrix}\,\,\,\,\h{for}\,\,\,\,\{j,k\}:=\Supp\,\aaa\,\,\,\,\h{with}\,\,\,\,j<k.
\leqno{(\rm A.2.4)}$$
In the case $\,[\aaa]\in\La_2$, $\,[\bb]\in\La_3\,$, we have
$$\aligned\dett(\aaa,\bb)&\in\{0\}\sqcup\{\pm\1 1\,{\pm}\,i\}\sqcup\{\pm\1 2,\,\,\pm\1 2i\},\\|\dett(\aaa,\bb)|&\in\{0,\,\sqrt{2},\,\,2\}.
\endaligned
\leqno{(\rm A.2.5)}$$
Note that $|\dett(\aaa,\bb)|$ is well-defined for $[\aaa],[\bb]$.
\ms
(v) For $i\in\{1,2,3\}$, $k\in\{2,3,6\}$, set
$$\aligned\La_j^m(\aaa)&:=\bl\{[\bb]\in\La_j\mid[\aaa]\ne[\bb],\,\,\mult(\aaa,\bb)=m\br\}\\\La^m(\aaa)&:=\h{$\bigsqcup$}_{j=1}^3\,\La_j^m(\aaa),\\ \La_j(\aaa)&:=\h{$\bigsqcup$}_m\,\La_j^m(\aaa).\endaligned
\leqno{(\rm A.2.6)}$$
For $[\aaa]\in\La_2$, set
$$\aligned'\!\La_2^6(\aaa)&:=\bl\{[\bb]\in\La_2^6\mid\Supp\,\aaa=\Supp\,\bb\br\},\\''\!\La_2^6(\aaa)&:=\bl\{[\bb]\in\La_2^6\mid\Supp\,\aaa\cup\Supp\,\bb=\{1,\dots,4\}\br\}.\endaligned
\leqno{(\rm A.2.7)}$$
For $[\aaa]\in\La_3$, set
$$\La_3^{m,d}(\aaa):=\bl\{[\bb]\in\La_3^m\mid{\rm diff}([\aaa],[\bb])=d\br\}.
\leqno{(\rm A.2.8)}$$
\ms
(vi) Using the action of $\Gamma$ in (A.1.4), we can define
$$\la_{j,j'}^m:=|\La_{j'}^m(\aaa)|,\q\la_j^m:=|\La^m(\aaa)|\,\bl(=\msum_{j'}\,\la_{j,j'}^m\br),\q\la_{j,j'}:=|\La_{j'}(\aaa)|\,\bl(=\msum_m\,\la_{j,j'}^m\br),
\leqno{(\rm A.2.9)}$$
independently of $[\aaa]\in\La_j$. We can similarly define
$${}'\!\la_{2,2}^6:=|{}'\!\La_2^6(\aaa)|,\q{}''\!\la_{2,2}^6:=|{}''\!\La_2^6(\aaa)|,\q\la_{3,3}^{m,d}:=|\La_3^{m,d}(\aaa)|,
\leqno{(\rm A.2.10)}$$
independently of $[\aaa]\in\La_2$ or $\La_3$.
\msn
{\bf A.3.~Classification of intersections.} For $[\aaa]\in\La_j$, $[\bb]\in\La_{j'}$ with $[\aaa]\ne[\bb]$, we can describe $\La_{j''}(\aaa,\bb)$ ($j''\in[1,3]$) together with $\la_{j,j'}^m$ for each $(j,j')\in[1,3]{\times}[1,3]$ as follows (where only the information of $\la_{j,j'}^m$ is noted in the case $j>j'$).
\ms
\vbox{\nin Case~(1,1) (that is, $(j,j')=(1,1)$).
\skn\q\q\,\,\,$\mult(\aaa,\bb)=6$, $\,|\La_2(\aaa,\bb)|=4$, $\,\,\,\La_2(\aaa,\bb)=\bl\{[\cc]\in\La_2\mid\Supp\,\cc=\Supp\,\aaa\sqcup\Supp\,\bb\br\}$.
\skn\q\q\,\,\,$\la_{1,1}^6=\tbinom{3}{1}=3$.}
\ms
\vbox{\nin Case~(1,2).
\ms
(a) If $\,\Supp\,\aaa\subset\Supp\,\bb,\,$ then $\,\,\mult(\aaa,\bb)=6$, $\,|\La_1(\aaa,\bb)|=1$, $\,|\La_2(\aaa,\bb)|=3$,
\sk
\q\q\q and $\,\lah(\aaa,\bb)\,$ is as in Case~(1,1) with $\aaa$, $\bb$ replaced appropriately.
\skn\q\q\,\,\,$\la_{1,2}^6=\tbinom{3}{1}\cdot 4=12$.}
\ms
\vbox{(b) If $\,\Supp\,\aaa\not\subset\Supp\,\bb,\,$ then $\,\,\mult(\aaa,\bb)=2\,$ (that is, $\La(\aaa,\bb)=\emptyset$).
\skn\q\q\,\,\,$\la_{1,2}^2=\tbinom{3}{2}\cdot 4=12$.}
\ms
\vbox{\nin Case~(1,3).
\skn\q\q\,\,\,$\mult(\aaa,\bb)=3$, $\,\,\La_3(\aaa,\bb)=\{[\cc]\}\,\,\,$ with $\,\,\,\pi^{(k)}[\bb]=\pi^{(k)}[\cc],\,\,\,\,\{k\}=\Supp\,\aaa$.
\skn\q\q\,\,\,$\la_{1,3}^3=|\La_3|=32$.}
\ms
\vbox{\nin Case~(2,1).
\skn\q\q\,\,\,$\la_{2,1}^2=\la_{2,1}^6=\tbinom{2}{1}=2$.}
\ms
\vbox{\nin Case~(2,2).
\ms
(a) If $\,|\Supp\,\aaa\cup\Supp\,\bb|=2,\,$ then $\,\mult(\aaa,\bb)=6$, $|\La_1(\aaa,\bb)|=2$, $|\La_2(\aaa,\bb)|=2$,
\sk
\q\q\q and $\,\lah(\aaa,\bb)\,$ is as in Case~(1,1) with $\aaa$, $\bb$ replaced appropriately.
\skn\q\q\,\,\,$'\la_{2,2}^6=4-1=3$.}
\ms
\vbox{(b) If $\,|\Supp\,\aaa\cup\Supp\,\bb|=3,\,$ then $\,\mult(\aaa,\bb)=3$, $\,\,\La_2(\aaa,\bb)=\{[\cc]\}$
\sk
\q\q\q with $\,\,\Supp\,\cc=(\Supp\,\aaa\cup\Supp\,\bb)\setminus(\Supp\,\aaa\cap\Supp\,\bb)$.
\skn\q\q\,\,\,$\la_{2,2}^3=\tbinom{2}{1}\cdot\tbinom{2}{1}\cdot 4=16$.}
\ms
\vbox{(c) If $\,|\Supp\,\aaa\cup\Supp\,\bb|=4$, $\,\prd(\aaa{+}\bb)=\pm\1 1,\,$ then $\,\mult(\aaa,\bb)=6$, $|\La_3(\aaa,\bb)|=4$
\sk
\q\q\q with $\,\,\La_3(\aaa,\bb)=\bl\{[\1\xi\1\aaa{+}\bb]\in\La_3\mid\xi\in\mmu_4\br\}$.
\skn\q\q\,\,\,$''\la_{2,2}^6=4/2=2$.}
\ms
\vbox{(d) If $\,|\Supp\,\aaa\cup\Supp\,\bb|=4$, $\,\prd(\aaa{+}\bb)=\pm\1 i,\,$ then $\,\mult(\aaa,\bb)=2$.
\skn\q\q\,\,\,$\la_{2,2}^2=4/2=2$.}
\ms
\vbox{\nin Case~(2,3). We may assume $a_{k_0}=b_{k_0}=1$ for some $k_0\in[1,4]$.
\ms
(a) If $\,\dett(\aaa,\bb)=0,\,$ then $\,\mult(\aaa,\bb)=6$, $|\La_2(\aaa,\bb)|=1$, $|\La_3(\aaa,\bb)|=3$,
\sk
\q\q\q and $\,\lah(\aaa,\bb)\,$ is as in Case~(2,2)(c) with $\aaa$, $\bb$ replaced appropriately.
\skn\q\q\,\,\,$\la_{2,3}^6=4\cdot 4/2=8$.}
\ms
\vbox{(b) If $\,|\dett(\aaa,\bb)|=\sqrt{2},\,$ then $\,\mult(\aaa,\bb)=3$, $\,\,\La_3(\aaa,\bb)=\{[\cc]\}\,$ with $\,{\rm rank}(\aaa,\bb,\cc)=2$
\skn\q\q\,\,\,$\la_{2,3}^3=2\cdot 4\cdot 4/2=16$.}
\ms
\vbox{(c) If $\,|\dett(\aaa,\bb)|=2,\,$ then $\,\mult(\aaa,\bb)=2$.
\skn\q\q\,\,\,$\la_{2,3}^2=4\cdot 4/2=8$.}
\ms
\vbox{\nin Case~(3,1).
\skn\q\q\,\,\,$\la_{3,1}^3=\tbinom{4}{1}=4$.}
\ms
\vbox{\nin Case~(3,2).
\skn\q\q\,\,\,$\la_{3,2}^2=\tbinom{4}{2}=6$, $\,\,\,\la_{3,2}^3=2\cdot\tbinom{4}{2}=12$, $\,\,\,\la_{3,2}^6=\tbinom{4}{2}=6$.}
\ms
\vbox{\nin Case~(3,3). We may assume $a_{k_0}=b_{k_0}=1$ for some $k_0\notin\Supp(\aaa-\bb)$.
\ms
(a) If $\,{\rm diff}([\aaa],[\bb])={\rm diff}(\aaa,\bb)=1,\,$ then $\,\mult(\aaa,\bb)=3$, $\,|\La_1(\aaa,\bb)|=1$,
\sk
\q\q\q and $\,\lah(\aaa,\bb)\,$ is as in Case~(1,3) with $\aaa$, $\bb$ replaced appropriately.
\skn\q\q\,\,\,$\la_{3,3}^{3,1}=\tbinom{4}{1}=4$.}
\ms
\vbox{(b) If $\,{\rm diff}([\aaa],[\bb])={\rm diff}(\aaa,\bb)=2$, $\,{\rm det}(\aaa-\bb,\bb)\ne 0,\,$ then $\,\mult(\aaa,\bb)=3$, $|\La_2(\aaa,\bb)|=1,$
\sk
\q\q\q and $\,\lah(\aaa,\bb)\,$ is as in Case~(2,3)(b) with $\aaa$, $\bb$ replaced appropriately.
\skn\q\q\,\,\,$\la_{3,3}^{3,2}=\tbinom{4}{2}\cdot 2=12$.}
\ms
\vbox{(c) If $\,{\rm diff}([\aaa],[\bb])={\rm diff}(\aaa,\bb)=2$, $\,{\rm det}(\aaa-\bb,\bb)=0,\,$ then $\,\mult(\aaa,\bb)=6$, $\,|\La_2(\aaa,\bb)|=2$,
\sk
\q\q\q$|\La_3(\aaa,\bb)|=2,\,$ and $\,\lah(\aaa,\bb)\,$ is as in Case~(2,2)(c) with $\aaa$, $\bb$ replaced appropriately.
\skn\q\q\,\,\,$\la_{3,3}^6=\tbinom{4}{2}\cdot 3\cdot\tfrac{1}{2}=9$.}
\ms
\vbox{(d) If $\,{\rm diff}([\aaa],[\bb])={\rm diff}(\aaa,\bb)=3,\,$ then $\,\mult(\aaa,\bb)=2$.
\skn\q\q\,\,\,$\la_{3,3}^2=3!=6$.}
\msn
{\bf Remark.} For Cases (2,3)(b) and (c), we use the following.
\msn
For $(a,b)\in\mmu_4{\times}\mmu_4$ and $e\in\{\pm\1 1\pm\1 i\}$, there are exactly two $(c,d)\in\mmu_4{\times}\mmu_4$ satisfying
$$\dett\begin{pmatrix}a&b\\ c&d\end{pmatrix}=e.$$
For $e\in\{\pm\1 2,\,\pm\1 2i\}$, there is only one $(c,d)\in\mmu_4{\times}\mmu_4$ satisfying the above relation.
\msn
{\bf A.4.~Conclusion.} Summarizing the above calculations, the table of the $\la_{j,j'}^m$ is as follows:
$$\begin{tabular}{|c|ccc|ccc|ccc|}
\hline
$\,j\,$ & & 1 & & & 2 & & & 3 &\\
\hline
$j'$ & 1 & 2 & 3 & 1 & 2 & 3 & 1 & 2 & 3\\
\hline
$m\,{=}\,2$ & 0 & 12& 0 & 2 & 2 & 8 & 0 & 6 & 6\\
$m\,{=}\,3$ & 0 & 0 & 32& 0 & 16& 16& 4 & 12& 16\\
$m\,{=}\,6$ & 3 & 12& 0 & 2 & 5 & 8 & 0 & 6 & 9\\
\hline\end{tabular}$$
\skn
Set as in (A.2.9)
$$\la_{j,j'}:=\msum_m\,\la_{j,j'}^m\,,\q\la_j^m:=\msum_{j'}\,\la_{j,j'}^m\,.$$
We see that the following relations hold for any $j\in[1,3]$:
$$\la_{j,1}=4-\delta_{j,1},\q\la_{j,2}=24-\delta_{j,2},\q\la_{j,3}=32-\delta_{j,3},
\leqno{(\rm A.4.1)}$$
$$\la_j^2=12,\q\la_j^3=2\cdot16,\q\la_j^6=5\cdot 3,
\leqno{(\rm A.4.2)}$$
(where $\delta_{j,k}=1$ if $j=k$, and $0$ otherwise) so that
$$\msum_{j'}\,\la_{j,j'}=\msum_m\,\la_j^m=60-1.$$
The above {\it double relations} are quite important to assure that the above calculations are correctly done and moreover we have counted {\it all\1} the singular points of $X_H$, where $X_H$ is a (reduced) projective line arrangement obtained by a general hyperplane section of $X$ in $\PP^3$. Note that (A.4.1) is closely related to (A.1.3), and (A.4.2) implies the independence of $\la^m_j/(m{-}1)$ on $j$, hence the latter coincides with the number of multiplicity $m$ points of $X_H$ contained in any fixed line in $X_H$, see also the proof of \cite[Proposition 5.1]{BDY} (and \cite{OT}).
\msn
{\bf Remarks.} (i) As in Case~(3,3)(d), the $\,[\bb]\in\La_3^2(\aaa)\,$ for $\,[\aaa]\in\La_3$ are given by
$$\bb=\ee\cdot\aaa\q\h{for}\q\ee=(e_k)\in(\mmu_4)^4\q\h{with}\q\{e_1,\dots,e_4\}=\mmu_4.$$
This implies that $\La_3$ has two (2)-connected components. Here a subset $J\subset\La$ is called $(m)$-{\it connected\1} if
$$X_H^J:=\mcup_{[\aaa]\in J}\,X_{H,\aaa}\,\,\,\,\h{with}\,\,\,\,X_{H,\al}:=X_H\cap X_{\al}$$
is connected via points of $X_H$ with multiplicity $m$ (that is, it is connected after deleting the singular points of $X_H$ with multiplicities different from $m$). Indeed, let $J'$ be a (2)-connected component of $\La_3$. The above assertion implies that $J'$ is stable by operations $[\aaa]\mapsto[\bb]$ defined by $\bb=\aaa\cdot\ee$ with $\ee=(e_k)$ satisfying
$$\{k\in[1,4]\mid e_k=1\}=2,\q\prd(\ee)=1.$$
(For instance, consider the composition of $(1,i,-1,-i)$ and the inverse of $(i,1,-1,-i)$, which produces $(-i,i,1,1)$.) Since $\prd(\aaa)$ remains invariant by these operations, we see that there are two (2)-connected components depending on $\prd(\aaa)=1$ or $-1$.
\ms
(ii) As in Case~(3,3)(a), the $\,[\bb]\in\La_3^{3,1}(\aaa)\,$ for $\,[\aaa]\in\La_3$ are given by
$$\bb=\ee\cdot\aaa\q\h{for}\q\ee=(e_k)\in(\mmu_4)^4\q\h{with}$$
$$\#\{k\mid e_k=1\}=3,\q\prd(\ee)=-1.$$
\ms
(iii) As in Case~(3,3)(b), the $\,[\bb]\in\La_3^{3,2}(\aaa)\,$ for $\,[\aaa]\in\La_3$ are given by
$$\bb=\ee\cdot\aaa\q\h{for}\q\ee=(e_k)\in(\mmu_4)^4\q\h{with}$$
$$\#\{k\mid e_k=1\}=2,\q\#\{k\mid e_k=\pm\1i\}=2,\q\prd(\ee)=1.$$
(This does not contradict Remark~(i) above, which treats ``connected components".)
\ms
(iv) As in Case~(3,3)(c), the $\,[\bb]\in\La_3^6(\aaa)\,$ for $\,[\aaa]\in\La_3$ are given by
$$\bb=\ee\cdot\aaa\q\h{for}\q\ee=(e_k)\in(\mmu_4)^4\q\h{with}$$
$$\#\{k\mid e_k=1\}=2,\q e_j=e_k\,\,\,\h{for}\,\,\,e_j\ne 1,\,\,e_k\ne 1.$$
\ms
(v) As in Case~(2,3)(a), the $\,[\bb]\in\La_3^6(\aaa)\,$ for $\,[\aaa]\in\La_2$ are given by the condition
$$\dett(\aaa,\bb)=0.$$
\ms
(vi) As in Case~(2,2)(a) and (c), the $\,[\bb]\in\La_2^6(\aaa)\,$ for $\,[\aaa]\in\La_2$ are given by the condition
$$|\Supp\,\aaa\cup\Supp\,\bb|=2\q\h{or}\q\prd(\aaa{+}\bb)=\pm\1 1.$$
Note that the last condition implies that $|\Supp\,\aaa\cup\Supp\,\bb|=4$.
\msn
{\bf A.5~Application.} We can apply Theorem~1 to the calculation of $H^1(F_{\!f},\C)_{\la}$ as in (2.5) for $\la=\exp(-2\pi i/6)$ in the case of $G_{31}$. Here we apply it to a general hyperplane section $X_H$ of $G_{31}$ as in (A.4). The subset $I\subset \La$ with $|I|=60/6=10$ is given, for instance, by
$$\aligned&(0,0,0,1),\\(-i,0,0,1),\q&(0,-i,0,1),\q(0,0,-i,1),\\(1,-1,0,0),\q&(1,0,-1,0),\q(0,1,-1,0),\\(i,i,1,1),\q&(i,1,i,1),\q(1,i,i,1).\endaligned$$
We can verify that the assumptions of Theorem~1 are satisfied by using the calculations in (A.3--4). Here $X_d$ is given by the first member, that is, $(0,0,0,1)$. Its restriction to $X_H$ intersects the lines defined by the next 3 members of $I$ at different points of $X_H$ with multiplicity 6, but there are no other intersections between the lines defined by these 10 members at multiplicity 6 points of $X_H$, see Remarks~(iv-vi) after (A.4) and the table in (A.4). So condition~(a) is satisfied.
\sk
We see moreover that $I^b$ is the union of $\La_1$ and the subset of $\La_2$ consisting of $[\aaa]\in\La_2$ with $\Supp\,\aaa\ni 4$. (Hence $|I^b|=4+3\cdot 4=16$.) Condition~(b) then follows from Remark~(i) after (A.4). Indeed, $\la_{2,3}^2\ne 0$ by the table in (A.4), and the two (2)-connected components of $\La_3$ in Remark~(i) after (A.4) are (2)-connected via an intermediate element of $\La_2$ as in Case~(2,3)(c); for instance, apply the latter to the case $\aaa=(1,-1,0,0)$ and $\bb=(1,1,1,1)$ or $(i,i,1,1)$, where $|\dett(\aaa,\bb)|=2$. (Note that a double point of $X_H^{I^c}$ which is also a double point of $X_H$ cannot be contained in $X_H^{I^b}$, and $\al_p\notin\Z$ unless $p$ is a point of $X_H$ with multiplicity 6.) Conditions~(c) and (d) follow from (A.4.2) and a remark after it. (Note that $\la^2_j+\la^3_j/2=28>|I^b|=16$ and $I^c=\{1,\dots,d\}\setminus I^b$.) So the hypotheses of Theorem~1 are verified.
\msn
{\bf Remarks.} (i) We have $H^1(F_{\!f},\C)_{\la}=0$ for $\la=\exp(-2\pi i/6)$ in this case. Indeed, it is rather easy to show the vanishing of the first cohomology of the Aomoto complex by using Remark~(ii) after (2.3), see also \cite[Lemma 3.1]{LY}, \cite[Section 1.5]{BDS}, \cite{BDY}, \cite{DS}.
\ms
(ii) We can show the above assertion also for $\la=-1$ by taking partitions $\La_j=\La'_j\sqcup\La''_j$ such that $|\La'_j|=|\La''_j|$ ($j\in[1,3]$). For instance, set for $[\aaa]\in\La_1$, $[\bb]\in\La_2$, $[\cc]\in\La_3$
$$\aligned{[\aaa]}\in\La'_1&\iff\Supp\,\aaa=\{3\}\,\,\,\h{or}\,\,\,\{4\},\\{[\bb]}\in\La'_2&\iff\h{$[\bb]\,$ is real},\,4\in\Supp\,\bb\,\,\,\,\,\h{or}\,\,\,\,\,\h{$[\bb]\,$ is not real},\,4\notin\Supp\,\bb,\\{[\cc]}\in\La'_3&\iff\prd(\cc)=1.\endaligned$$
Here $[\bb]$ is called {\it real\1} if it is represented by an element of $\R^4$. Setting $I:=\bigsqcup_j\La'_j$, we have $|I^b|=6$, where $X_d$ is given by $[(0,0,0,1)]$. Indeed, we have $\al_p=0$ for any singular point $p$ of $X_H$ with multiplicity 6 except for the intersection point $p_0$ of the lines defined by $(0,0,0,1)$ and $(0,0,1,0)$ where $\al_{p_0}=1$.
\ms
(iii) It seems rather difficult to find the subset $I$ for $\la=\exp(\pm 2\pi i/3)$, since there are so many triple points of $X_H$, see Cases (1,3), (2,2)(b), (2,3)(b), (3,3)(a), (3,3)(b) as well as Remarks~(ii) and (iii) after (A.4). It seems still difficult even if we replace $P_1$ with $P_2$, where $P$ is the pole order filtration. There is also a certain difficulty in extending Theorem~1 to the case for $P_2$.


\begin{thebibliography}{MPP}
\bibitem[BDY]{BDY} Bailet, P., Dimca, A.\ and Yoshinaga, M., A vanishing result for the first twisted cohomology of affine varieties and applications to line arrangements (arXiv:1705.06022), Manuscripta Math.\ (2018)
\bibitem[BBD]{BBD} Beilinson, A., Bernstein, J.\ and Deligne, P., Faisceaux pervers, Ast\'erisque 100, Soc.\ Math.\ France, Paris, 1982.
\bibitem[Br]{Br} Brieskorn, E., Sur les groupes de tresses, S\'eminaire Bourbaki, 24\`eme ann\'ee (1971/1972), Exp. No. 401, Lect.\ Notes in Math.\
317, Springer, Berlin, 1973, pp. 21--44.
\bibitem[BS]{BS} Budur, N.\ and Saito, M., Jumping coefficients and spectrum of a hyperplane arrangement, Math.\ Ann.\ 347 (2010), 545--579. 
\bibitem[BDS]{BDS} Budur, N., Dimca, A.\ and Saito, M., First Milnor cohomology of hyperplane arrangements, Contemp.\ Math.\ 538, AMS, Providence, RI, 2011, pp.~279--292.
\bibitem[BSY]{BSY} Budur, N., Saito, M.\ and Yuzvinsky, S., On the local zeta functions and the b-functions of certain hyperplane arrangements (With an appendix by Willem Veys), J.\ London Math.\ Soc.\ (2) 84 (2011), 631--648. 
\bibitem[De]{De} Deligne, P., Th\'eorie de Hodge II, Publ.\ Math.\ IHES, 40 (1971), 5--58.
\bibitem[Di1]{Di1} Dimca, A., Singularities and Topology of Hypersurfaces, Springer, Berlin, 1992.
\bibitem[Di2]{Di2} Dimca, A., Sheaves in topology, Springer, Berlin, 2004.
\bibitem[Di3]{Di3} Dimca, A., On the Milnor monodromy of the irreducible complex reflection arrangements (arXiv:1606.04048), J.\ Inst.\ Math.\ Jussieu (2017)
\bibitem[DS]{DS} Dimca, A.\ and Sticlaru, G., On the Milnor monodromy of the exceptional reflection arrangement of type $G_{31}$, Documenta Math.\ 23 (2018), 1--14.
\bibitem[ESV]{ESV} Esnault, H., Schechtman, V.\ and Viehweg, E., Cohomology of local systems on the complement of hyperplanes, Inv.\ Math.\ 109 (1992), 557--561.
\bibitem[Gr]{Gr} Grothendieck, A., El\'ements de g\'eom\'etrie alg\'ebrique III-1, Publ.\ Math.\ IHES 11 (1961).
\bibitem[HR]{HR} Hoge, T.\ and R\"ohrle, G., On supersolvable reflection arrangements, Proc.\ AMS, 142 (2014), 3787--3799.
\bibitem[JR]{JR} Jozsa, R.\ and Rice, J., On the cohomology ring of hyperplane complements, Proc.\ AMS, 113 (1991), 973--981.
\bibitem[LY]{LY} Libgober, A. and Yuzvinsky, S., Cohomology of the Orlik-Solomon algebras and local systems, Compos.\ Math.\ 121 (2000), 337--361.
\bibitem[MPP]{MPP} Ma\u cinic, A., Papadima, S.\ and Popescu, C.R., Modular equalities for complex reflexion arrangements, Documenta Math.\ 22 (2017), 135--150.
\bibitem[OS]{OS} Orlik, P.\ and Solomon, L., Combinatorics and topology of complements of hyperplanes, Inv.\ Math.\ 56 (1980), 167--189.
\bibitem[OT]{OT} Orlik, P.\ and Terao, H., Arrangements of Hyperplanes, Springer, Berlin, 1992.
\bibitem[SaK]{SaK} Saito, K., Theory of logarithmic differential forms and logarithmic vector fields, J.\ Fac.\ Sci.\ Univ.\ Tokyo Sect.\ IA Math., 27(2) (1980), 265--291.
\bibitem[Sa1]{mhp} Saito, M., Modules de Hodge polarisables, Publ. RIMS, Kyoto Univ. 24 (1988), 849--995.
\bibitem[Sa2]{mhm} Saito, M., Mixed Hodge modules, Publ. RIMS, Kyoto Univ.\ 26 (1990), 221--333.
\bibitem[Sa3]{bCM} Saito, M., Multiplier ideals, $b$-function, and spectrum of a hypersurface singularity, Compos.\ Math.\ 143 (2007), 1050--1068.
\bibitem[Sa4]{bha} Saito, M., Bernstein-Sato polynomials of hyperplane arrangements, Selecta Math.\ (N.S.) 22 (2016), 2017--2057.
\bibitem[Sa5]{ex} Saito, M., Hilbert series of graded Milnor algebras and roots of Bernstein-Sato polynomials, arXiv:1509.06288.
\bibitem[STV]{STV} Schechtman, V., Terao, H. and Varchenko, A., Local systems over complements of hyperplanes and the Kac-Kazhdan conditions for singular vectors, J. Pure Appl.\ Algebra 100 (1995), 93--102.
\bibitem[Wo]{Wo} Wotzlaw, L., Intersection cohomology of hypersurfaces, Ph.D. dissertation, Humboldt Universit\"at zu Berlin, August 2006.
\end{thebibliography}
\end{document}